\RequirePackage{fix-cm}
\documentclass[smallextended]{svjour3}       
\smartqed  

\usepackage{marvosym}

\renewcommand{\thesection}{\arabic{section}}

\usepackage{color}
\usepackage{xcolor}

\usepackage{graphicx}
\begin{document}

\title{Tempered Fractional Multistable  Motion and Tempered Multifractional Stable  Motion 
}

\titlerunning{Fractional multistable  motion and multifractional stable motion}        

\author{Xiequan Fan   \and
 Jacques L\'{e}vy V\'{e}hel
}

\authorrunning{X. Fan and  J. L\'{e}vy V\'{e}hel } 

\institute{X. Fan  (\Letter)   \at
Center for Applied Mathematics,
Tianjin University, Tianjin,  China \\
             \email{fanxiequan@hotmail.com}
           \and
           J. L\'{e}vy V\'{e}hel  (\Letter) \at
            Regularity Team, Inria, France\\
              \email{jacques.levy-vehel@inria.fr}
}

\date{Received: date / Accepted: date}

\maketitle

\begin{abstract} This work defines two classes of processes, that we term
{\it tempered fractional multistable motion} and {\it tempered multifractional stable motion}. They are extensions of fractional multistable  motion  and multifractional stable  motion, respectively,  obtained by adding
 an exponential tempering to the integrands. We investigate certain basic features of these  processes, including scaling property, tail probabilities,  absolute  moment, sample path properties, pointwise H\"{o}lder exponent, H\"{o}lder continuity  of quasi norm, (strong) localisability and semi-long-range dependence structure. These processes may provide useful models for data that exhibit both dependence and varying local regularity/intensity of jumps.

\keywords{Stable processes \and Multistable processes  \and Multifractional processes  \and Sample paths\and Long-range dependence\and Localisability }

\subclass{ 60G52 \and 60G18 \and 60G22 \and  60G17 \and 60E07  }
\end{abstract}

\section{Introduction}
Linear fractional stable  motion (LFSM) can be  represented by the stochastic integral of a symmetric $\alpha$-stable random measure $dZ_{\alpha}(x)$, that is
\begin{eqnarray}\label{ssfl2}
  X(t) =\int_{-\infty}^{\infty}\Big[ (t-x)_+^{H -\frac1{\alpha }} - (-x)_+^{H-\frac1{\alpha }}\Big]dZ_{\alpha}(x), \ \ \ \ \ t \in \mathbf{R},
\end{eqnarray}
where $0<\alpha \leq2, 0< H <1,$   $(x)_+=\max\{x, 0\}$ and $0^0=0$.
See for example \mbox{Samorodnitsky} and Taqqu \cite{ST94}. This stochastic process   has two important features. It is
self-similar with Hurst parameter $H,$ i.e.\ for any $c>0,$ $t_1,...,t_d \in  \mathbf{R},$
$$\big(X(c\, t_1),..., X(c\, t_d) \big)\stackrel{d}{=}\big( c^HX(t_1),..., c^HX(t_d) \big),$$ and it has stationary increments, i.e.,\ for any $\tau \in \mathbf{R},$
$\big(X(t)-X(0), \ -\infty < t < \infty \big )\stackrel{d}{=}\big( X(\tau+t)-X(\tau), \ -\infty < t < \infty \big),$
where $\stackrel{d}{=} $ indicates equality in distribution.
Because its
increments can exhibit the heavy-tailed analog of long-range dependence (see Watkins et al.\ \cite{WCH05}),
the model is useful in practice to model, for example, financial data, internet traffic, noise on telephone line, signal processing and atomospheric noise,
see Nolan \cite{N10} for many references.

There exist at least three extensions of  LFSM, i.e.,\ linear multifractional stable
motion (LmFSM),  linear fractional multistable
motion (LFmSM) and  linear tempered  fractional stable motion (LTFSM).
Stoev and Taqqu \cite{ST04,ST05} first introduced LmFSM by replacing the self-similarity parameter $H$ in the integral representation of the LFSM  by a time-varying function $H_t$.
Stoev and Taqqu have examined the effect of the regularity of the function $H_t$ on the local structure of the process. They also showed that under certain H\"{o}lder regularity conditions on the function $H_t$, the LmFSM is locally equivalent to a LFSM, in the sense of finite-dimensional distributions.
Thus LmFSM is  a locally self-similar stochastic process. Whereas the LFSM   is always continuous in probability, this is not in general the case for  LmFSM. Stoev and Taqqu have obtained necessary and sufficient conditions for the continuity in probability of the LmFSM.
Falconer and L\'{e}vy V\'{e}hel \cite{FL09}   defined the second  model extension of LFSM, called LFmSM. LFmSM  behaves locally like  linear fractional $\alpha(t)$-stable  motion close to time $t,$ in the sense that the local scaling limits are   linear fractional $\alpha(t)$-stable  motions, but where the stability index $\alpha(t)$
varies with $t$.  This extension allows one to account for the fact that the nature of irregularity, including the stability level, may vary in time.
See also Falconer  and  Liu \cite{FL12}  where the $\alpha$-stable random measure in (\ref{ssfl2}) has been  replaced by a time-varying $\alpha(t)$-multistable random measure.
Recently, \mbox{Meerschaert} and Sabzikar \cite{MS16}  defined the third extension, termed LTFSM, by  adding an exponential tempering to the
power-law kernel in a LFSM. They showed that the   LTFSM exhibits semi-long-range dependence,   and therefore
provides a useful alternative model for data that
exhibit strong dependence.

In view of trying to combine the properties of both LFmSM and LTFSM, we define in this work a new stochastic process by adding an exponentiel tempering to the power-law kernel of LFmSM. Our {\it linear tempered fractional multistable motion} (LTFmSM) is thus
  an extension of  LFmSM and LTFSM. In particular,
 linear tempered fractional multistable  motion  behaves locally like the linear fractional $\alpha(t)$-stable  motion with
stability index $\alpha(t)$ that varies in time $t,$ and  it exhibits semi-long-range dependence structure as LTFSM does.
Similarly,  to combine the properties of both LmFSM and LTFSM, we define another new stochastic process, called  {\it linear tempered multifractional stable motion} (LTmFSM), by adding an exponentiel tempering to the power-law kernel of LmFSM.  This new process is also of semi-long-range dependence structure.
We investigate basic properties of the two new processes, including  scaling properties, tail probabilities,  absolute  moment,
sample path properties, pointwise H\"{o}lder exponent, H\"{o}lder continuity of quasi norm and (strong) localisability.
Such properties  are important and have been widely studied. For instance,
Falconer and  Liu \cite{FL12} have  investigated  sample path properties, localisability and  strong  localisability of LFmSM;
Le Gu\'{e}vel and L\'{e}vy V\'{e}hel \cite{GLL13a} have investigated the pointwise H\"{o}lder exponent of LFmSM; Ayache and Hamonier \cite{AH14} have examined  the fine path properties of LmFSM;
Meerschaert  and Sabzikar \cite{MS16} have studied scaling properties, sample path properties and H\"{o}lder continuity of quasi norm of LTFSM.

The reader will note that, in this work, our emphasis is on the properties that set apart LTFmSM and LTmFSM, rather than
on their common ones. Further work is needed to introduce and study linear tempered multifractional multistable motion (LTmFmSM). We believe that studying the specific properties of LTFmSM and LTmFSM will be helpful for future investigation of LTmFmSM.

The remainder of this paper is organized as follows.
In Section \ref{sec2}, we define the linear tempered fractional multistable motion  and the linear tempered multifractional stable motion.
In Section \ref{sec3},   we elucidate the dependence structure of the two stochastic processes.
In Sections \ref{sec4} - \ref{sec8}, we analyze  their properties.

\section{Definitions of LTFmSM and LTmFSM}\label{sec2}
Throughout this paper, for given $0<a \leq  b \leq2,$ the function $\alpha: \mathbf{R} \longrightarrow [a, b]$
will be a Lebesgue measurable function that will play the role of a varying  stability index.  We recall the definition of
variable exponent  Lebesgue space:
$$\mathcal{F}_\alpha:=\{f: f\ \textrm{is measurable with}\ ||f||_\alpha < \infty\}$$
where
\begin{eqnarray}\label{seminorm}
||f||_\alpha:=\Bigg\{\lambda>0: \ \int_{-\infty}^\infty\Big|\frac{f(x)}{\lambda} \Big|^{\alpha(x)}dx =1   \Bigg\} .
\end{eqnarray}
Then $||\cdot||_\alpha$ is a quasinorm.

 Falconer  and  Liu \cite{FL12} defined the \textit{multistable stochastic integral}
$I(f):=\int f(x) dM_{\alpha}(x), f \in \mathcal{F}_\alpha,$ by specifying the finite-dimensional distribution of $I.$
Here and after, $dM_{\alpha}(x)$ stands for the multistable measure, which is an independently scattered symmetric random measure.
Assume  $\alpha(x) \in [a, b] \subset (0, 2]$. Given $f_1, f_2,...,f_d \in \mathcal{F}_\alpha,$ Falconer  and  Liu defined a probability distribution
on the vector $\big(I(f_1), I(f_2),...,I(f_d)\big)\in \mathbf{R}^d$ by the following characteristic function
\begin{eqnarray*}
 \mathbf{E}\Big[e^{i\sum_{k=1}^d \theta_k I(f_k)  } \Big]  =\  \exp\Bigg\{ - \int_{-\infty}^{\infty} \Big|\sum_{k=1}^n\theta_k  f_k(x)   \Big|^{\alpha(x)} dx\Bigg\}.
\end{eqnarray*}
The essential point is that $\alpha(x)$ may vary with $x.$ With the definition  of multistable stochastic integral, Falconer  and  Liu \cite{FL12}
(cf.\ Proposition 4.3 therein) defined linear fractional multistable
motion (LFmSM)
\begin{eqnarray} \label{ed2sd}
X(t)=\int_{-\infty}^{\infty}\Big[  (t-x)_+^{H -\frac1{\alpha(x)}} -  (-x)_+^{H-\frac1{\alpha(x)}}\Big]dM_{\alpha}(x).
\end{eqnarray}
They also investigated some basic properties of LFmSM, such as localisability and strong localisability.

By adding  an exponential tempering to the power-law kernel in  LFSM (\ref{ssfl2}), that is
\begin{eqnarray} \label{ed2d}
  X_{H,\alpha,\lambda}(t):=\int_{-\infty}^{\infty}\Big[ e^{-\lambda (t-x)_+}(t-x)_+^{H -\frac1{\alpha }} - e^{-\lambda (-x)_+}(-x)_+^{H-\frac1{\alpha }}\Big]dZ_{\alpha}(x),
\end{eqnarray}
 $\lambda>0, 0<\alpha < 2$ and $    H > 0,$
Meerschaert and Sabzikar \cite{MS16} recently defined the so-called linear tempered  fractional stable motion
(LTFSM).  They showed that LTFSM  is short memory, but its increments behave like long memory when $\lambda$ is very small.
 Thus LTFSM exhibits semi-long-range dependence structure, and   it
provides a useful alternative model for data that
exhibit strong dependence.

Similarly, by adding  an exponential tempering to the power-law kernel in a LFmSM (\ref{ed2sd}), we  define
the following linear tempered  fractional multistable motion. Such process
 is an extension of both LFmSM and LTFSM mentioned above.

\begin{definition} Let $\alpha(x)  \in [a, b] \subset (0, 2]$ be a continuous function  on $\mathbf{R}.$
 Given an independently scattered symmetric  multistable random measure $dM_{\alpha}(x)$ on $\mathbf{R}$, the multistable stochastic integral
\begin{eqnarray}\label{defi2.1}
  X_{H,\alpha(x),\lambda}(t):=\int_{-\infty}^{\infty}\Big[ e^{-\lambda (t-x)_+}(t-x)_+^{H -\frac1{\alpha(x)}} - e^{-\lambda (-x)_+}(-x)_+^{H-\frac1{\alpha(x)}}\Big]dM_{\alpha}(x)
\end{eqnarray}
with   $0<H <1, \lambda >0, (x)_+=\max\{x, 0\},$ and $0^\gamma=0, \gamma \in \mathbf{R},$ will be called a
\emph{linear tempered  fractional multistable motion} (LTFmSM).
\end{definition}
\begin{remark}
  With the exponential tempering, we can also define  \emph{multistable Yaglom noise}
$$Y_{H,\alpha(x),\lambda}(t)= \int_{-\infty}^{\infty}\Big[ e^{-\lambda(t-x)_+}(t-x)_+^{H-\frac1{\alpha(x)}}  \Big] d M_{\alpha}(x), \ \ \ \ \lambda>0. $$
In particular, when $ \alpha (x) \equiv 1/H \in (0, 2],$ multistable Yaglom noise is known as Ornstein-Uhlenbeck process, see Example 3.6.3 of Samorodnitsky and Taqqu \cite{ST94}. When $\alpha (x)\equiv \alpha$ for some constant $\alpha,$ multistable Yaglom noise is called stable Yaglom noise, see Meerschaert  and Sabzikar \cite{MS16}.  It is obvious that  fractional multistable Yaglom noise
is  a multistable stochastic integral. It is also easy to see that
$$X_{H,\alpha(x),\lambda}(t) =Y_{H,\alpha(x),\lambda}(t) -Y_{H,\alpha(x),\lambda}(0) , \ \ \ \ \ \lambda > 0.$$
\end{remark}

Denote by
\begin{eqnarray}\label{funcgs}
G_{H ,\alpha(x),\lambda }(t, x)= e^{-\lambda (t-x)_+}(t-x)_+^{H -\frac1{\alpha(x)}} - e^{-\lambda (-x)_+}(-x)_+^{H -\frac1{\alpha(x)}}, \ \ \ \ \ \lambda > 0.
\end{eqnarray}
It is easy to check that the function $G_{H ,\alpha(x),\lambda }(t, x)$ belong to $\mathcal{F}_\alpha, $ so that
LTFmSM is well defined. Moreover, by the definition of multistable integral (cf.\ Falconer  and  Liu \cite{FL12}),  the characteristic function of  $X_{H,\alpha(x),\lambda}(t)$
is given as follows:
\begin{eqnarray}\label{defnl}
 \mathbf{E}\Big[e^{i\sum_{k=1}^d \theta_k X_{H,\alpha(x),\lambda}(t_k)  } \Big]  =\  \exp\Bigg\{ - \int_{-\infty}^{\infty} \Big|\sum_{k=1}^d\theta_k  G_{H ,\alpha(x),\lambda }(t_k, x)   \Big|^{\alpha(x)} dx\Bigg\}.
\end{eqnarray}

Similarly, when the Hurst parameter $H$ of (\ref{ed2d}) varies with time $t$, we have another extension of  LTFSM.
\begin{definition} Let   $H_t \in [a, b]  $  be a continuous function  on $\mathbf{R}.$
 Given an independent scattered S$\alpha$S stable random measure $dZ_{\alpha }(x)$ on $\mathbf{R}$
with  control measure $dx$,  the  stable stochastic integral
\begin{eqnarray}\label{defin01}
X_{H_t,\alpha ,\lambda }(t):=\int_{-\infty}^{\infty}\Big[ e^{-\lambda(t-x)_+}(t-x)_+^{H_t-\frac1{\alpha}} - e^{-\lambda(-x)_+}(-x)_+^{H_t-\frac1{\alpha }}\Big] d Z_{\alpha }( x)
\end{eqnarray}
with $0<\alpha  \leq 2,$  $\lambda > 0, (x)_+=\max\{x, 0\},$ and  $0^\gamma=0, \gamma \in \mathbf{R},$  will be called a
\emph{linear tempered multifractional  stable motion} (LTmFSM).
\end{definition}

Denote
 $$G_{H_t ,\alpha,\lambda }(t, x)= e^{-\lambda (t-x)_+}(t-x)_+^{H_t -\frac1{\alpha}} - e^{-\lambda (-x)_+}(-x)_+^{H_t -\frac1{\alpha}}, \ \ \ \ \ \lambda \geq 0. $$
By the definition of stable integral (cf.\ Samorodnitsky  and Taqqu \cite{ST94}),  the characteristic function of  $X_{H_t,\alpha,\lambda}(t)$
is given as follows:
\begin{eqnarray*}
 \mathbf{E}\Big[e^{i\sum_{k=1}^n \theta_k X_{H_t,\alpha,\lambda}(t_k)  } \Big]  =\  \exp\Bigg\{ - \int_{-\infty}^{\infty} \Big|\sum_{k=1}^n\theta_k   G_{H_t ,\alpha,\lambda }(t_k, x)   \Big|^{\alpha} dx\Bigg\}.
\end{eqnarray*}
The characteristic function of  $X_{H_t,\alpha,\lambda}(t)$
is given as follows:
\begin{eqnarray*}
 \mathbf{E}\Big[e^{i\sum_{k=1}^n \theta_k X_{H_t,\alpha,\lambda}(t_k)  } \Big]  =\  \exp\Bigg\{ - \int_{-\infty}^{\infty} \Big|\sum_{k=1}^n\theta_k   G_{H_t ,\alpha,\lambda }(t_k, x)   \Big|^{\alpha} dx\Bigg\}.
\end{eqnarray*}

\section{Dependence structure of LTFmSM and LTmFSM}\label{sec3}

In this section, we study the behaviour of increments of LTFmSM and LTmFSM, usually termed the ``noise'' of these processes.

Denote by
\[
 Y (t) =X (t+1) -X (t)  \ \ \  \textrm{for integers}   \  -\infty<t < \infty
\]
the noise of the  processes $X$.
Astrauskas et al.\ \cite{ALT91} studied the dependence structure of linear fractional stable motion using the following nonparametric
measure of dependence (see also Meerschaert and Sabzikar \cite{MS16}).
Define
\[
 R_{t_1}(t)=R(\theta_1, \theta_2,  t_1, t_1  +t ):= \mathbf{E}\Big[ e^{i(\theta_1Y(t_1) + \theta_2Y(t_1+t))}    \Big]    - \mathbf{E}\Big[ e^{i \theta_1Y(t_1)}    \Big] \mathbf{E}\Big[ e^{i  \theta_2Y(t_1+t) }  \Big]
\]
for $t_1, t , \theta_1,  \theta_2 \in \mathbf{R}.$ If we also define
\begin{eqnarray}
  I(\theta_1, \theta_2, t_1, t_1+t) &=& \log\Big( \mathbf{E}\Big[ e^{i \theta_1Y(t_1)}    \Big]\Big)  + \log\Big(  \mathbf{E}\Big[ e^{i  \theta_2Y(t_1+t) }  \Big] \Big) \nonumber \\
  &&  - \log\Big( \mathbf{E}\Big[ e^{i(\theta_1Y(t_1) + \theta_2Y(t_1+t))}    \Big] \Big),\nonumber
 \end{eqnarray}
then we have
\begin{equation}
 R_{t_1}(t) = K(\theta_1, \theta_2, t_1, t_1+t)\Big( e^{-I(\theta_1, \theta_2, t_1, t_1+t)} -1\Big), \label{rtdef}
 \end{equation}
where
$$K(\theta_1, \theta_2, t_1, t_1+t ) =\mathbf{E}\Big[ e^{i \theta_1Y(t_1)}    \Big] \mathbf{E}\Big[ e^{i  \theta_2Y(t_1+t) }  \Big].$$
In particular, for stationary processes, $ R_{t_1}(t)$ does not depend on $t_1,$ see  Meerschaert and Sabzikar \cite{MS16}.
In this case, we denote $R_{t_1}(t)$ by $R(t)$ for simplicity. Note however that the increments of the two processes that we define in this work
are not stationary in general.

We first recall the dependence structure of LTFSM.
 Given two real-valued functions $f(t), g(t)$ on $\mathbf{\mathbf{R}}$,
we will write $$f(t)\preceq g(t)$$ if $|f(t)/g(t)| \leq C_1$ for all $t>0$
sufficiently large and some $0< C_1 <\infty.$
In particular, if $f(t)\preceq g(t)$ and $g(t)\preceq f(t),$
we  will write $$f(t) \asymp g(t).$$
Thus  $f(t) \asymp g(t)$ is equivalent to   $C_1 \leq|f(t)/ g(t)|  \leq C_2 $ for all $t>0$
sufficiently large and some $0< C_1< C_2 <\infty.$
With these notations,
Meerschaert and Sabzikar \cite{MS16} recently proved that if $ \lambda >0$  and $0< \alpha \leq1,$    then  TFSN has the following property
\[
R(t) \asymp  e^{-\lambda \alpha t }  t^{H\alpha -1}
\]
for  $\theta_1 \theta_2 \neq 0.$ Meerschaert and Sabzikar \cite{MS16} also proved that if $ \lambda >0, 1<  \alpha <2$  and $ \frac1\alpha < H,$    then TFSN   has the following property
\[
R(t) \asymp  e^{-\lambda   t }  t^{H-\frac 1\alpha }
\]
for  $\theta_1 \theta_2 \neq 0.$

\subsection{Dependence structure of LTFmSM}

 The following two theorems show that  LTFmSM  and LTFSM share  the similar  dependence structure.

\begin{definition} Given an LTFmSM defined by (\ref{defi2.1}), we define the tempered fractional multistable noise (TFmSN)
\begin{equation} \label{fsfsfb}
 Y_{H, \alpha(x), \lambda}(t):=X_{H, \alpha(x), \lambda}(t+1) -X_{H, \alpha(x), \lambda}(t)
\end{equation}
for integers $-\infty<t < \infty.$
\end{definition}

In particular, if $\alpha(x)\equiv \alpha$ for a constant $\alpha \in (0, 2]$, then the TFmSN reduces to  the   tempered fractional  stable noise, see Meerschaert and Sabzikar \cite{MS16}.

\begin{proposition}\label{th31}
Let $\alpha(x)  \in [a, b] \subset (0, 1)$ be a continuous function  on $\mathbf{R}.$
Let $Y_{H, \alpha(x), \lambda}(t)$ be the tempered fractional multistable noise (\ref{fsfsfb}). Recall $R_{t_1}(t)$ defined by (\ref{rtdef})
with $Y(t)= Y_{H, \alpha(x), \lambda}(t)$. Assume  $ \lambda >0.$
Then
\begin{equation}\label{sfddf}
e^{-\lambda b t }  t^{H a  -1}  \preceq R_{t_1}(t) \preceq  e^{-\lambda a t }  t^{Hb -1}
\end{equation}
for any $t_1 \in \mathbf{R}$ and $\theta_1 \theta_2 \neq 0.$
\end{proposition}
\emph{Proof.} By the definition (\ref{defi2.1}), TFmSN has  the following representation
\[
Y_{H, \alpha(x), \lambda}(t)= \int_{-\infty}^{\infty}\Big[e^{-\lambda (t+1-x)_+}(t+1-x)_+^{H -\frac1{\alpha(x)}}  - e^{-\lambda (t-x)_+}(t-x)_+^{H -\frac1{\alpha(x)}}  \Big]dM_{\alpha}(x).
\]
Define $g_t(x)=e^{-\lambda (t-x)_+}(t-x)_+^{H-\frac1{\alpha(x)}}$ for $t\in \mathbf{R}$ and write
\begin{eqnarray}
I(\theta_1, \theta_2, t_1, t_1+t)&=& \int_{-\infty}^{\infty}\Big| \theta_1 [g_{t_1+1}(x) -g_{t_1}(x)] + \theta_2  [g_{t_1+t+1}(x) -g_{t_1+t}(x)] \Big|^{\alpha(x)}   dx \nonumber \\
 && \ \ - \int_{-\infty}^{\infty}\Big| \theta_1 [g_{t_1+1}(x) -g_{t_1}(x)]\Big|^{\alpha(x)}   dx  \nonumber \\
 && \ \ \ \ \   -  \int_{-\infty}^{\infty}\Big|   \theta_2  [g_{t_1+t+1}(x) -g_{t_1+t}(x)] \Big|^{\alpha(x)}   dx  \nonumber \\
 &=& I_1(t) + I_2( t), \label{f26iz}
\end{eqnarray}
where
\begin{eqnarray*}
I_1( t )&=& \int_{-\infty}^{t_1} \bigg( \Big| \theta_1 [g_{t_1+1}(x) -g_{t_1}(x)] + \theta_2  [g_{t_1+t+1}(x) -g_{t_1+t}(x)] \Big|^{\alpha(x)}    \nonumber \\
 &&\ \ \ \ \ \ \ \ \ \ \ \ -  \Big| \theta_1 [g_{t_1+1}(x) -g_{t_1}(x)]\Big|^{\alpha(x)}
 -   \Big|   \theta_2  [g_{t_1+t+1}(x) -g_{t_1+t}(x)] \Big|^{\alpha(x)}   \bigg) dx
\end{eqnarray*}
and
\begin{eqnarray*}
I_2(  t)&=& \int_{t_1}^{t_1+1} \bigg( \Big| \theta_1 g_{t_1+1}(x) + \theta_2  [g_{t_1+t+1}(x) -g_{t_1+t}(x)] \Big|^{\alpha(x)}    \nonumber \\
 &&\ \ \ \ \ \ \ \ \ \ \ \ \   -  \Big| \theta_1 g_{t_1+1}(x) \Big|^{\alpha(x)}
 -   \Big|   \theta_2  [g_{t_1+t+1}(x) -g_{t_1+t}(x)] \Big|^{\alpha(x)}  \bigg) dx .
\end{eqnarray*}
Using the following inequalities
\begin{eqnarray} \label{ghllgm}
0 \leq |x_1|^\alpha + |x_2|^\alpha  - |x_1+x_2|^\alpha   \leq  2\, |x_2|^\alpha
\end{eqnarray}
for all $x_1, x_2 \in \mathbf{R}$ and $0<\alpha \leq 1,$   we obtain
\begin{eqnarray}\label{gmsnom}
 I_1( t ) \leq 0 \ \ \ \ \  \textrm{and} \ \ \  \  \ \  I_2( t ) \leq 0.
\end{eqnarray}
First, we give an estimation for $I_1( t )$.  By (\ref{ghllgm}),  it is easy to see that  for $t \geq 1,$
\begin{eqnarray*}
 |I_1( t )|  &\leq&  2  \int_{-\infty}^{t_1}\Big|   \theta_2  [g_{t_1+t+1}(x) -g_{t_1+t}(x)]   \Big|^{\alpha(x)}   dx \\
&\leq&  2 \Big(|\theta_2|^a + |\theta_2|^b \Big)e^{-\lambda a t }  t^{Hb -1} \int_{-\infty}^{t_1}\Big|   [g_{t_1+t+1}(x) -g_{t_1+t}(x)]e^{ \lambda   t }  t^{\frac{1}{\alpha(x)}  -H } \Big|^{\alpha(x)}   dx.
\end{eqnarray*}
Notice that $H\alpha(x)\leq1.$
For $x \leq t_1$ and $t>1$, we deduce that
\begin{eqnarray}
&&\Big|   [g_{t_1+t+1}(x) -g_{t_1+t}(x)]e^{ \lambda   t }  t^{\frac{1}{\alpha(x)}  -H } \Big|^{\alpha(x)}  \nonumber \\
&&= \bigg|e^{-\lambda   (t_1-x)} \bigg(e^{-\lambda  } \Big(1+ \frac{1+ t_1-x}{t}\Big)^{H-\frac 1{\alpha(x)}}  - \Big(1+\frac{ t_1-x}{t}\Big)^{H-\frac 1{\alpha(x)}} \bigg)  \bigg|^{\alpha(x)} \nonumber\\
&&\leq   e^{-\lambda a (t_1-x)} (1 + e^{-\lambda  })^b    \bigg( \Big(1+ \frac{1+ t_1-x}{t}\Big)^{H\alpha(x)- 1 }  + \Big(1+\frac{ t_1-x}{t}\Big)^{H\alpha(x)- 1 } \bigg)   \nonumber \\
&&\leq  2 e^{-\lambda a (t_1-x)} (1 + e^{-\lambda  })^b   \nonumber\\
&&=: F_\lambda(x) .  \label{sfdsg}
\end{eqnarray}
Thus
\begin{eqnarray*}
|I_1( t )| \, e^{ \lambda a t }  t^{1-Hb} &\leq&  2 (|\theta_2|^a + |\theta_2|^b )  \int_{-\infty}^{t_1} F_\lambda(x)   dx \nonumber \\
&\leq&  C_1 (|\theta_2|^a + |\theta_2|^b ) ,
\end{eqnarray*}
 where $C_1>0$ depends only on $a, b$ and $\lambda.$
Hence
\begin{eqnarray}
  |I_1( t )|
&\leq&  C_1 (|\theta_2|^a + |\theta_2|^b )  e^{- \lambda a t }  t^{Hb -1} .  \label{ineq021}
\end{eqnarray}
Next for $I_2(  t)$, we have the following estimation.
Using inequality (\ref{ghllgm}) again, we obtain
\begin{eqnarray}\label{ghlm}
| I_2( t)|    \leq  2\int_{t_1}^{t_1+1}\Big|   \theta_2  [g_{t_1+t+1}(x) -g_{t_1+t}(x)]  \Big|^{\alpha(x)}   dx.
\end{eqnarray}
Applying the mean value theorem to see that for   $t\geq2$  and any $x \in (t_1, t_1+1)$, we have
\begin{eqnarray}
&& \Big| g_{t_1+t+1}(x) -g_{t_1+t}(x) \Big| \nonumber \\
&& \leq \Big|- \lambda e^{-\lambda(u-x)}(u-x)^{H-\frac1{\alpha(x)}} +(H-\frac1{\alpha(x)})e^{-\lambda(u-x)}(u-x)^{H-\frac1{\alpha(x)}-1} \Big|  \nonumber \\
&&\leq  e^{-\lambda(t-1)}\Big( \lambda(t-1)^{H-\frac1{\alpha(x)}} +(\frac1{\alpha(x)} -H)(t-1)^{H-\frac1{\alpha(x)}-1} \Big) \nonumber \\
&&\leq  e^{-\lambda(t-1)}  (\frac1{\alpha(x)} -H+ \lambda) (t-1)^{H-\frac1{\alpha(x)}} ,\nonumber
\end{eqnarray}
where $u \in (t_1+t, t_1+t+1).$
Returning to (\ref{ghlm}), we get
\begin{eqnarray}
|I_2(  t)|    &\leq& 2(|\theta_2|^a + |\theta_2|^b ) \int_{t_1}^{t_1+1}\Big| e^{-\lambda(t-1)}  (\frac1{\alpha(x)} -H+ \lambda) |t-1|^{H-\frac1{\alpha(x)}}  \Big|^{\alpha(x)}   dx \nonumber  \\
&\leq& C_2 (|\theta_2|^a + |\theta_2|^b )  e^{-\lambda a t }  t^{Hb -1}  \label{ineq022}
\end{eqnarray}
for large $t,$ where $C_2>0$ depends only on $a, b, H$ and $\lambda.$
 Combining  the inequalities (\ref{f26iz}), (\ref{gmsnom}),(\ref{ineq021})  and (\ref{ineq022}) together, we obtain
\begin{eqnarray}\label{inep24}
0\leq - I(\theta_1, \theta_2, t_1, t_1+t) \leq C_3   e^{-\lambda a t }  t^{Hb -1} \ \ \ \ \textrm{as}   \ t\rightarrow \infty,
\end{eqnarray}
where $C_3$ does not depend  on $t.$
Using the following  equality
\begin{eqnarray*}
  |x_1|^\alpha + |x_2|^\alpha  - |x_1+x_2|^\alpha   =  |x_2|^\alpha -\frac{\alpha}{|x_1 + \theta  x_2|^{1-\alpha}} |x_2|
\end{eqnarray*}
for all $x_1, x_2 \neq 0 $ with $|x_2|$ small enough,  $0<\alpha < 1$  and some $|\theta|\leq 1,$ we obtain for any $x_1 \neq 0$ and $0<\alpha < 1,$
\begin{eqnarray}\label{ghm2}
 |x_1|^\alpha + |x_2|^\alpha  - |x_1+x_2|^\alpha   \sim  |x_2|^\alpha
\end{eqnarray}
as $  x_2 \rightarrow 0.$
It is easy to see that
for $  t_1 \leq x \leq t_1+1$ and $t>2$,
\begin{eqnarray}
&& \lim_{t\rightarrow \infty}     [g_{t_1+t+1}(x) -g_{t_1+t}(x)]e^{ \lambda   t }  t^{\frac{1}{\alpha(x)}  -H }   \nonumber \\
&&= \lim_{t\rightarrow \infty} e^{-\lambda   (t_1-x)} \bigg(e^{-\lambda  } \Big(1+ \frac{1+ t_1-x}{t}\Big)^{H-\frac 1{\alpha(x)}}  - \Big(1+\frac{ t_1-x}{t}\Big)_+^{H-\frac 1{\alpha(x)}} \bigg)    \nonumber\\
&&=  e^{-\lambda   (t_1-x)} \big(e^{-\lambda  } - 1\big)       \label{fmgdgfs}
\end{eqnarray}
and
\begin{eqnarray*}
   \Big|   [g_{t_1+t+1}(x) -g_{t_1+t}(x)]    \Big|^{\alpha(x)}
  \leq  \Big|e^{-\lambda   (t_1-x)} \big(2+ e^{-\lambda  } \big)  \Big|^{\alpha(x)} e^{ - a \lambda   t }  t^{1- H a}.
\end{eqnarray*}
Thus $    [g_{t_1+t+1}(x) -g_{t_1+t}(x)],$ $  t_1 \leq x \leq t_1+1,$ converges uniformly to $0$ as $t\rightarrow \infty.$
  Applying the dominated convergence theorem yields for  $\theta_1 \theta_2 \neq 0,$ we have
\begin{eqnarray}
| I_2( t)| \succeq \int^{t_1+1 }_{t_1}\Big|   \theta_2  [g_{t_1+t+1}(x) -g_{t_1+t}(x)]   \Big|^{\alpha(x)}   dx
 \end{eqnarray}
and
\begin{eqnarray}
\liminf_{t\rightarrow \infty} |I_2( t )| \, e^{ \lambda b t }  t^{1-Ha} &\geq& \lim_{t\rightarrow \infty}\int_{t_1}^{t_1+1}  \Big|  \theta_2 [g_{t_1+t+1}(x) -g_{t_1+t}(x)]e^{ \lambda   t }  t^{\frac{1}{\alpha(x)}  -H } \Big|^{\alpha(x)}  dx \nonumber \\
&=& \int_{t_1}^{t_1+1}\Big|   \theta_2  \big[e^{-\lambda   (t_1-x)} \big(1- e^{-\lambda  } \big) \big]   \Big|^{\alpha(x)}   dx . \label{fmdggdg}
\end{eqnarray}
Then (\ref{f26iz}),  (\ref{gmsnom}) and (\ref{fmdggdg}) implies that for all large $t$,
\begin{eqnarray} \label{inepfs25}
- I(\theta_1, \theta_2, t_1, t_1+t) \geq -I_2( t ) =| I_2( t)|  \geq \frac12 C_2 \,  e^{ -\lambda b t }  t^{Ha -1},
\end{eqnarray}
where $C_2=\int_{t_1}^{t_1+1}\big|   \theta_2  \big[e^{-\lambda   (t_1-x)} \big(1- e^{-\lambda  } \big) \big]   \big|^{\alpha(x)}   dx$ does not depend  on $t.$ Combining (\ref{inep24}) and (\ref{inepfs25}) together, we have
 \begin{eqnarray} \label{f26fds}
 e^{ -\lambda b t }  t^{Ha -1} \preceq  I(\theta_1, \theta_2, t_1, t_1+t)  \preceq    e^{ -\lambda a t }  t^{H b -1}
\end{eqnarray}
for   $\theta_1 \theta_2 \neq 0.$
It is easy to see that
\begin{eqnarray}
 &&  K(\theta_1, \theta_2, t_1, t_1+t) \nonumber  \\
 &&= \exp\bigg\{ - \int_{-\infty}^{1}\Big|   \theta_1  [e^{-\lambda (1-u)_+}(1-u)_+^{H -\frac1{\alpha(t_1+u)}}  - e^{-\lambda (-u)_+}(-u)_+^{H -\frac1{\alpha(t_1+u)}} ] \Big|^{\alpha(t_1+u)}   du\bigg\} \nonumber  \\
 && \times \exp\bigg\{   -\int_{-\infty}^{1}\Big|   \theta_2  [e^{-\lambda ( 1-u)_+}(1-u)_+^{H -\frac1{\alpha(t_1+t+u)}}  - e^{-\lambda (-u)_+}(-u)_+^{H -\frac1{\alpha(t_1+t+u)}} ]  \Big|^{\alpha(t_1+t+u)}   du \bigg\}  \nonumber \\
 && \geq  \exp\bigg\{ - 2 \big( |\theta_1|^a + |\theta_2|^b \big)\int_{-\infty}^{1} M_\lambda(u)    du \bigg\}, \nonumber
 \end{eqnarray}
where $$M_\lambda(u)=   e^{-\lambda a(1-u)_+} \Big((1-u)_+^{Ha -1} + (1-u)_+^{Hb -1} \Big) + e^{-\lambda a (-u)_+}  \Big(( -u)_+^{Ha -1} + ( -u)_+^{Hb -1} \Big)  $$
 is integrable on  $ (-\infty, 1]$ with respect to $u$,  and that
 $$ K(\theta_1, \theta_2, t_1, t_1+t)  \leq 1.$$
Since $I(\theta_1, \theta_2, t_1, t_1+t)   \rightarrow 0 $ as $t\rightarrow \infty, $
it follows from (\ref{rtdef}) that $R_{t_1}(t) \sim -K(\theta_1, \theta_2, t_1, t_1+t)I(\theta_1, \theta_2, t_1, t_1+t).  $
Hence (\ref{sfddf}) follows by (\ref{f26fds}).    \qed

\vspace{0.2cm}

\begin{proposition}\label{th32}
Let $\alpha(x)  \in [a, b] \subset (1, 2]$ be a continuous function  on $\mathbf{R}.$
Let $Y_{H, \alpha(x), \lambda}(t)$ be the tempered fractional multistable noise (\ref{fsfsfb}). Recall $R_{t_1}(t)$ defined by (\ref{rtdef})
with $Y(t)=Y_{H, \alpha(x), \lambda}(t)$. Assume $ \lambda >0$
and $1/a < H <1.$
Then
\begin{equation}\label{ineq39}
e^{-\lambda  t }t^{H-\frac1 a } \preceq R_{t_1}(t) \preceq e^{-\lambda  t }t^{H-\frac1 b }
\end{equation}
for any $t_1 \in \mathbf{R}$ and $\theta_1 \theta_2 \neq 0.$
\end{proposition}
\emph{Proof.} Recall $ I_1(t)$  and $ I_2( t) $ defined by (\ref{f26iz}). Notice that
\begin{eqnarray}\label{ghllfgm}
 |x_1+x_2|^\alpha -  |x_1|^\alpha  -  |x_2|^\alpha \sim \alpha \,  x_1  \,  |x_2|^{\alpha-1}
\end{eqnarray}
for all $ x_2 \neq 0, x_1\rightarrow 0$ and $1< \alpha \leq 2. $
First, we give an estimation for $ I_1(t ).$  It is easy to see that  for $ x\leq t_1+1,$ $1/a < H$ and large $t$,
\begin{eqnarray}
&&\Big|   [g_{t_1+t+1}(x) -g_{t_1+t}(x)]e^{ \lambda   t }  t^{\frac{1}{\alpha(x)}  -H } \Big|   \nonumber \\
&&= \bigg|e^{-\lambda   (t_1-x)} \bigg(e^{-\lambda  } \Big(1+ \frac{1+ t_1-x}{t}\Big)^{H-\frac 1{\alpha(x)}}  - \Big(1+\frac{ t_1-x}{t}\Big)^{H-\frac 1{\alpha(x)}} \bigg)  \bigg| \nonumber \\
&&\leq   e^{-\lambda  (t_1-x)} (1 + e^{-\lambda  })   \bigg( \Big(1+ \frac{1+ t_1-x}{t}\Big)^{H-\frac 1{\alpha(x)}}  + \Big(1+\frac{ t_1-x}{t}\Big)^{H-\frac 1{\alpha(x)}} \bigg)   \nonumber \\
&&\leq  2 e^{-\lambda   (t_1-x)} (1 + e^{-\lambda  }) (2+ t_1-x  )^{H-\frac 1{\alpha(x)}}  \label{fine524}
\end{eqnarray}
and
\begin{eqnarray}
&&    [g_{t_1+t+1}(x) -g_{t_1+t}(x)]e^{ \lambda   t }  t^{\frac{1}{\alpha(x)}  -H }    \nonumber \\
&&=  e^{-\lambda  (t_1-x)}    \bigg( \Big(e^{-\lambda/(H-\frac 1{\alpha(x)})} (1+\frac1 t) + e^{-\lambda/(H-\frac 1{\alpha(x)})} \frac{  t_1-x}{t}\Big)^{H-\frac 1{\alpha(x)}}  - \Big(1+\frac{ t_1-x}{t}\Big)^{H-\frac 1{\alpha(x)}} \bigg)   \nonumber \\
&&\leq  e^{-\lambda  (t_1-x)}    \bigg( \Big(1+ e^{-\lambda/(H-\frac 1{\alpha(x)})} \frac{  t_1-x}{t}\Big)^{H-\frac 1{\alpha(x)}}  - \Big(1+\frac{ t_1-x}{t}\Big)^{H-\frac 1{\alpha(x)}} \bigg)   \nonumber \\
&&\leq 0. \label{sgjmf}
\end{eqnarray}
 Then (\ref{fine524}) and (\ref{sgjmf}) together implies that  for $ x\leq t_1+1,$ $1/a < H$ and large $t$,
\begin{eqnarray}
0&\geq&g_{t_1+t+1}(x) -g_{t_1+t}(x)     \nonumber \\
&\geq&  -2 (1 + e^{-\lambda  }) e^{-\lambda t}t^{H- \frac1{b}} e^{-\lambda   (t_1-x)} (2+ t_1-x  )^{H-\frac 1{b}}\  \nonumber \\
  &\geq&  -2 (1 + e^{-\lambda  }) e^{2\lambda + \frac1 b -H   } \Big( \frac{H-\frac1 b }{b} \Big)^{H-\frac 1{b}}\ e^{-\lambda t}t^{H- \frac1{b}} . \nonumber
\end{eqnarray}
Thus $[g_{t_1+t+1}(x) -g_{t_1+t}(x)]$ is negative and converges to $0$  uniformly for $x \in (-\infty, t_1]$ as  $t \rightarrow \infty.$
By (\ref{ghllfgm}), we obtain for large $t$,
\begin{eqnarray*}
&&  |I_1( t)|  \\
&&\leq   2 \int_{-\infty}^{t_1}\alpha(x)\Big|   \theta_2  [g_{t_1+t+1}(x) -g_{t_1+t}(x)]   \Big| \Big| \theta_1 [g_{t_1+1}(x) -g_{t_1}(x)] \Big|^{\alpha(x)-1}   dx \\
&&\leq 4| \theta_2|\max\Big\{ |\theta_1|^{a-1}, |\theta_1|^{b-1} \Big \}   \int_{-\infty}^{t_1}\Big| g_{t_1+t+1}(x) -g_{t_1+t}(x)  \Big|   \Big|  g_{t_1+1}(x) -g_{t_1}(x)  \Big|^{\alpha(x)-1}  dx.
\end{eqnarray*}
Therefore, for  large $t$ and $x\leq t_1$,
\begin{eqnarray}
& &   |I_1( t)|  \nonumber\\
 &&\leq  4| \theta_2|\max\{ |\theta_1|^{a-1}, |\theta_1|^{b-1} \} \nonumber \\
 && \ \ \ \times\int_{-\infty}^{t_1} e^{ -\lambda   t }  t^{H-\frac1{\alpha(x)}} (1 +e^{-\lambda  }) \Big(2 +t_1-x\Big) ^{H-\frac{1}{\alpha(x)}  } e^{-\lambda  (t_1-x)}  \Big|  g_{t_1+1}(x) -g_{t_1}(x)  \Big|^{\alpha(x)-1}  dx \nonumber \\
 &&\leq  4| \theta_2|\max\{ |\theta_1|^{a-1}, |\theta_1|^{b-1} \} e^{ -\lambda   t }  t^{H-\frac1{b}}\nonumber \\
 && \ \ \ \times\int_{-\infty}^{t_1} (1 +e^{-\lambda  }) \Big(2 +t_1-x\Big) ^{H-\frac{1}{\alpha(x)}  } e^{-\lambda  (t_1-x)}  \Big|  g_{t_1+1}(x) -g_{t_1}(x)  \Big|^{\alpha(x)-1}  dx  .\label{inqgs1}
\end{eqnarray}
Recall $g_t(x)=e^{-\lambda (t-x)_+}(t-x)_+^{H-\frac1{\alpha(x)}},$ and that
\begin{eqnarray*}
 \Big| g_{t_1 +1}(x) -g_{t_1 }(x) \Big|^{\alpha(x)-1}
  \leq \Big| g_{t_1 +1}(x)\Big|^{\alpha(x)-1}+  \Big|g_{t_1 }(x) \Big|^{\alpha(x)-1}
\end{eqnarray*}
(cf.   (\ref{ghllgm}) for the last inequality). Since $\alpha(x)-1 \leq b-1 <1 $ and $H>\frac1{a}\geq\frac1{\alpha(x)},$ from (\ref{inqgs1}), we obtain
\begin{eqnarray}
   |I_1( t )|  \leq  C_1 | \theta_2|\max\Big\{ |\theta_1|^{a-1}, |\theta_1|^{b-1} \Big \}e^{ -\lambda   t }  t^{H-\frac1{b}} ,\label{ineq02g1}
\end{eqnarray}
where $C_1$ does not depend on $t.$
Next, we give an estimation for $I_2(  t).$ Using  (\ref{ghllfgm}) again,   we obtain for large $t,$
\begin{eqnarray*}
 |I_2(  t)| &=& \int_{t_1}^{t_1+1} \Bigg| \Big| \theta_1 g_{t_1+1}(x) + \theta_2  [g_{t_1+t+1}(x) -g_{t_1+t}(x)] \Big|^{\alpha(x)}    \nonumber \\
 &&\ \ \ \ \ \ \ \ \ \ \ \ \   -  \Big| \theta_1 g_{t_1+1}(x) \Big|^{\alpha(x)}
 -   \Big|   \theta_2  [g_{t_1+t+1}(x) -g_{t_1+t}(x)] \Big|^{\alpha(x)}  \Bigg| dx \\
 &\leq& 2\int_{t_1}^{t_1+1} \alpha(x)\Big|   \theta_2  [g_{t_1+t+1}(x) -g_{t_1+t}(x)] \Big|\Big| \theta_1 g_{t_1+1}(x) \Big|^{\alpha(x)-1}  dx \\
 &\leq& 4 | \theta_2|\max\Big\{ |\theta_1|^{a-1}, |\theta_1|^{b-1} \Big \} \int_{t_1}^{t_1+1}  \Big|  g_{t_1+t+1}(x) -g_{t_1+t}(x)  \Big|\Big|  g_{t_1+1}(x) \Big|^{\alpha(x)-1}  dx .
\end{eqnarray*}
By (\ref{fine524}),  it follows that for large $t,$
\begin{eqnarray}
 |I_2( t)|
 &\leq& 4 | \theta_2|\max\Big\{ |\theta_1|^{a-1}, |\theta_1|^{b-1} \Big \} e^{ -\lambda   t }  t^{H-\frac{1}{b} } \nonumber \\
 &&   \times \int_{t_1}^{t_1+1}  e^{-\lambda   (t_1-x)} (1 + e^{-\lambda  }) (2+ t_1-x  )^{H-\frac 1{\alpha(x)}}   \Big|  g_{t_1+1}(x) \Big|^{\alpha(x)-1}  dx \nonumber  \\
&\leq& C_2 | \theta_2|\max\Big\{ |\theta_1|^{a-1}, |\theta_1|^{b-1} \Big\}\  e^{-\lambda t}  t^{H-\frac1{b}} , \label{ineq02g2}
\end{eqnarray}
where $C_2$ does not depend  on $t.$
Therefore, from   (\ref{ineq02g1}) and (\ref{ineq02g2}), for large $t,$
\begin{eqnarray}\label{gknm25}
 |I(\theta_1, \theta_2, t_1, t_1+t)|  \leq C_3  e^{-\lambda t}  t^{H-\frac1{b}} .
\end{eqnarray}
where $C_3$ does not depend  on $t.$

By (\ref{ghllfgm}),   we have
\begin{eqnarray*}
   |I_2( t)| \succeq  \int_{t_1}^{t_1+1}\alpha(x)   \Big| \theta_2  [g_{t_1+t+1}(x) -g_{t_1+t}(x)] \Big|    \Big| \theta_1  g_{t_1+1}(x)   \Big|^{\alpha(x)-1}   dx .
\end{eqnarray*}
Applying  (\ref{fmgdgfs}) and the dominated convergence theorem yields
\begin{eqnarray*}
& &\liminf_{t\rightarrow \infty} | I_2( t)| e^{\lambda  t }t^{\frac1 a -H}   \nonumber\\
 &&\geq  \lim_{t\rightarrow \infty}  | \theta_2|\min\{ |\theta_1|^{a-1}, |\theta_1|^{b-1} \} \nonumber \\
 && \ \ \ \times\int_{t_1}^{t_1+1} \Big|   [g_{t_1+t+1}(x) -g_{t_1+t}(x)]e^{ \lambda   t }  t^{\frac{1}{\alpha(x)}  -H } \Big|  \Big|  g_{t_1+1}(x)     \Big|^{\alpha(x)-1}  dx \nonumber \\
 &&=  | \theta_2|\min\{ |\theta_1|^{a-1}, |\theta_1|^{b-1} \}\int_{t_1}^{t_1+1} e^{-\lambda   (t_1-x)} \big(1-e^{-\lambda  })  \Big|  g_{t_1+1}(x)    \Big|^{\alpha(x)-1}  dx  .
\end{eqnarray*}
Thus for $\theta_1 \theta_2 \neq 0,$
\begin{eqnarray}\label{fgjm5}
  |I_2( t)| \succeq  e^{-\lambda  t }t^{H-\frac1 a}.
\end{eqnarray}
Notice that (\ref{ghllfgm}) and (\ref{sgjmf}) implies that
\begin{eqnarray}\label{fghlmm5}
I_1(t)I_2(t) \geq 0
\end{eqnarray}
 for large $t.$
Combining (\ref{gknm25}) and (\ref{fgjm5}) together, we have for   $\theta_1 \theta_2 \neq 0,$
 \begin{eqnarray} \label{f26s}
e^{-\lambda  t }t^{H-\frac1 a} \preceq |I_2( t)|  \preceq | I(\theta_1, \theta_2, t_1, t_1+t) | \preceq   e^{-\lambda  t }t^{H-\frac1 b}.
\end{eqnarray}

Since $I(\theta_1, \theta_2, t_1, t_1+t)   \rightarrow 0 $ as $t\rightarrow \infty, $
it follows from (\ref{rtdef}) that $R_{t_1}(t) \sim -K(\theta_1, \theta_2, t_1, t_1+t)I(\theta_1, \theta_2, t_1, t_1+t);  $
hence (\ref{ineq39}) holds.   \qed

\subsection{Dependence structure of LTmFSM}
In this section, we consider the increment of  LTmFSM. The following two theorems extend the dependence structure of LTFSM  to the case of  LTmFSM.

\begin{definition} Given an  LTmFSM defined by (\ref{defin01}), we define the tempered multifractional stable noise (TmFSN)
\begin{equation}\label{ineq30}
 Y_{H_t, \alpha, \lambda}(t):=X_{H_{t+1}, \alpha, \lambda}(t+1) -X_{H_t, \alpha, \lambda}(t)
\end{equation}
for integers $-\infty<t < \infty.$
\end{definition}

In particular, if $H_t\equiv H$ for a constant $H \in (0, 1)$, then the TmFSN reduces to  the  tempered fractional  stable noise.
The next theorem shows that LTmFSM has a dependence structure more general than that of LTFSM.
\begin{proposition}\label{th41}
Let   $H_t \in [a, b]  $  be a continuous function  on $\mathbf{R}.$
Let $Y_{H_t, \alpha, \lambda}(t)$ be a tempered multifractional stable noise (\ref{ineq30}) for some $0< \alpha <1$.
Recall $R_{t_1}(t)$ defined by (\ref{rtdef}) with $Y(t)=Y_{H_t, \alpha, \lambda}(t)$.
Assume $ \lambda >0.$
Then
\begin{equation}\label{ineq31}
R_{t_1}(t) \asymp e^{- \lambda \alpha t }  t^{\alpha H_t -1}
\end{equation}
  for  $\theta_1 \theta_2 \neq 0.$
\end{proposition}
\emph{Proof.} By the definition (\ref{defin01}), TmFSN has  the following representation
\[
Y_{H_t, \alpha, \lambda}(t)= \int_{-\infty}^{\infty}\Big[e^{-\lambda (t+1-x)_+}(t+1-x)_+^{H_{t+1} -\frac1{\alpha}}  - e^{-\lambda (t-x)_+}(t-x)_+^{H_t -\frac1{\alpha}}  \Big] dZ_{\alpha}(x).
\]
Define $h_t(x)=(t-x)_+^{H_t-\frac1{\alpha}}e^{-\lambda (t-x)_+}$ for $t\in \mathbf{R}$ and write
\begin{eqnarray}
I(\theta_1, \theta_2, t_1, t_1+t)&=& \int_{-\infty}^{\infty}\Big| \theta_1 [h_{t_1+1}(x) -h_{t_1}(x)] + \theta_2  [h_{t_1+t+1}(x) -h_{t_1+t}(x)] \Big|^{\alpha}   dx \nonumber \\
 && \ \ - \int_{-\infty}^{\infty}\Big| \theta_1 [h_{t_1+1}(x) -h_{t_1}(x)]\Big|^{\alpha}   dx  \nonumber \\
 && \ \ \ \ \   -  \int_{-\infty}^{\infty}\Big|   \theta_2  [h_{t_1+t+1}(x) -h_{t_1+t}(x)] \Big|^{\alpha}   dx  \nonumber \\
 &=& I_3(t ) + I_4(  t),  \label{ineq33}
\end{eqnarray}
where
\begin{eqnarray*}
I_3( t )&=& \int_{-\infty}^{t_1} \bigg( \Big| \theta_1 [h_{t_1+1}(x) -h_{t_1}(x)] + \theta_2  [h_{t_1+t+1}(x) -h_{t_1+t}(x)] \Big|^{\alpha}     \nonumber \\
 &&\ \ \ \ \ \ \ \ \ \ \    -  \Big| \theta_1 [h_{t_1+1}(x) -h_{t_1}(x)]\Big|^{\alpha}
 -   \Big|   \theta_2  [h_{t_1+t+1}(x) -h_{t_1+t}(x)] \Big|^{\alpha} \bigg)  dx
\end{eqnarray*}
and
\begin{eqnarray*}
I_4( t)&=& \int_{t_1}^{t_1+1} \bigg( \Big| \theta_1 h_{t_1+1}(x) + \theta_2  [h_{t_1+t+1}(x) -h_{t_1+t}(x)] \Big|^{\alpha}    \nonumber \\
 &&\ \ \ \ \ \ \ \ \ \ \ \ \   -  \Big| \theta_1 h_{t_1+1}(x) \Big|^{\alpha}
 -   \Big|   \theta_2  [h_{t_1+t+1}(x) -h_{t_1+t}(x)] \Big|^{\alpha}  \bigg) dx .
\end{eqnarray*}
Using  (\ref{ghllgm}) again, we obtain
\begin{eqnarray}\label{fgh2h}
 I_3( t )  \leq 0 \ \ \ \ \ \textrm{and} \ \ \ \     I_4( t )  \leq 0  .
\end{eqnarray}
First, we give an estimation for $I_3( t)$.  For large $t,$
\begin{eqnarray*}
| I_3( t  )|  &\leq&  2  \int_{-\infty}^{t_1}\Big|   \theta_2  [ h_{t_1+t+1}(x) -h_{t_1+t}(x)]   \Big|^{\alpha}   dx \\
&\leq&  2 |\theta_2|^\alpha e^{-\lambda \alpha t }  t^{ \alpha H_t-1} \int_{-\infty}^{t_1}\Big|   [h_{t_1+t+1}(x) -h_{t_1+t}(x)]e^{ \lambda   t }  t^{\frac{1}{\alpha}  -H_t } \Big|^{\alpha}   dx.
\end{eqnarray*}
Recall $H_t \in [a, b] . $ It is easy to see that
for $x \leq t_1$ and $t>1$,
\begin{eqnarray*}
&&\Big|   [h_{t_1+t+1}(x) -h_{t_1+t}(x)]e^{ \lambda   t }  t^{\frac{1}{\alpha}  -H_t } \Big|^{\alpha}  \nonumber \\
&&= \bigg|e^{-\lambda   (t_1-x)} \bigg(e^{-\lambda  } \Big(1+ \frac{1+ t_1-x}{t}\Big)^{H_t-\frac 1{\alpha }}  - \Big(1+\frac{ t_1-x}{t}\Big)^{H_t-\frac 1{\alpha}} \bigg)  \bigg|^{\alpha} \\
&&\leq   e^{-\lambda \alpha (t_1-x)} (1 + e^{-\lambda  })^\alpha   \bigg( \Big(1+ \frac{1+ t_1-x}{t}\Big)^{H_t\alpha - 1 }  + \Big(1+\frac{ t_1-x}{t}\Big)^{H_t\alpha- 1 } \bigg)    \\
&&\leq  e^{-\lambda \alpha (t_1-x)} (1 + e^{-\lambda  })^\alpha \max\bigg\{2, \,   \big(2+ t_1-x   \big)^{b\alpha - 1 }  + \big(1+t_1-x  \big)^{b\alpha- 1 }  \bigg\}\\
&&:= F_\lambda(x).
\end{eqnarray*}
Thus
\begin{eqnarray*}
|I_3( t )|e^{ \lambda \alpha t }  t^{1-\alpha H_t} &\leq&  2 |\theta_2|^\alpha    \int_{-\infty}^{t_1} F_\lambda(x) dx \nonumber \\
&\leq&  C_1 \, |\theta_2|^\alpha   ,
\end{eqnarray*}
 where $C_1>0$ depends only on $ \alpha, b$ and $\lambda.$
Hence
\begin{eqnarray}\label{ineq35}
|I_3( t )|
&\leq&  C_1 \,  |\theta_2|^\alpha   e^{- \lambda \alpha t }  t^{\alpha H_t -1} .
\end{eqnarray}
Next for $I_4(  t)$, we have the following estimation.
Using inequality (\ref{ghllgm}) again, we obtain
\begin{eqnarray}\label{ineq36}
\Big|I_4(  t ) \Big|   \leq  2\int_{t_1}^{t_1+1}\Big|   \theta_2  [h_{t_1+t+1}(x) -h_{t_1+t}(x)]  \Big|^{\alpha}   dx.
\end{eqnarray}
Applying the mean value theorem to see that for   $t\geq2$  and any $x \in (t_1, t_1+1)$, we have
\begin{eqnarray}
&& \Big| h_{t_1+t+1}(x) -h_{t_1+t}(x) \Big| \nonumber \\
&& \leq \Big|- \lambda e^{-\lambda(u-x)}(u-x)^{H_t-\frac1{\alpha}} +(H_t-\frac1{\alpha})e^{-\lambda(u-x)}(u-x)^{H_t-\frac1{\alpha}-1} \Big|  \nonumber \\
&&\leq  e^{-\lambda(t-1)}\Big( \lambda(t-1)^{H_t-\frac1{\alpha}} +(\frac1{\alpha} -H_t)(t-1)^{H_t-\frac1{\alpha}-1} \Big) \nonumber \\
&&\leq  e^{-\lambda(t-1)}  (\frac1{\alpha} -H_t+ \lambda) (t-1)^{H_t-\frac1{\alpha}} ,\label{fnsggh}
\end{eqnarray}
where $u \in (t_1+t, t_1+t+1).$
Returning to (\ref{ineq36}), we get
\begin{eqnarray}
\Big|I_4(  t) \Big|   &\leq& 2 |\theta_2|^ \alpha   \int_{t_1}^{t_1+1}\Big| e^{-\lambda(t-1)}  (\frac1{\alpha} -H_t+ \lambda) |t-1|^{H_t-\frac1{\alpha}}  \Big|^{\alpha}   dx \nonumber  \\
&\leq& C_2 |\theta_2|^\alpha     e^{-\lambda \alpha t }  t^{\alpha H_t -1}   \label{ineq38}
\end{eqnarray}
for large $t,$ where $C_2>0$ depends only on $\alpha, b$ and $\lambda.$
 Combining  the inequalities (\ref{ineq33}), (\ref{fgh2h}), (\ref{ineq35}) and (\ref{ineq38}) together, we obtain
\begin{eqnarray}\label{fh963h}
0\leq -I(\theta_1, \theta_2, t_1, t_1+t)  \leq C_3  |\theta_2|^\alpha     e^{-\lambda \alpha t }  t^{\alpha H_t -1}
\end{eqnarray}
for large $t$, where $C_3$ does not depend  on $t.$
By (\ref{ghm2}), it holds for $t\rightarrow \infty,$
\[
| I_4( t)| \succeq \int_{t_1}^{t_1+1}\Big|   \theta_2  [h_{t_1+t+1}(x) -h_{t_1+t}(x)]   \Big|^{\alpha }   dx.
\]
Similar to (\ref{fmgdgfs}), it is easy to see that
for $t_1\leq x \leq t_1+1$ and $t>2$,
\begin{eqnarray}
  \lim_{t\rightarrow \infty}  \Big|   [h_{t_1+t+1}(x) -h_{t_1+t}(x)]e^{ \lambda   t }  t^{\frac{1}{\alpha}  -H_t } \Big|^{\alpha}  =  \big|e^{-\lambda   (t_1-x)} \big( 1 -e^{-\lambda  }   \big)  \big|^{\alpha}  .  \label{fmcwgdg}
\end{eqnarray}
 Applying the dominated convergence theorem yields
\begin{eqnarray}
\liminf_{t\rightarrow \infty} |I_4( t )| \, e^{ \lambda \alpha t }  t^{1-H_t \alpha } &\geq&
\lim_{t\rightarrow \infty}\int_{t_1}^{t_1+1}\Big|   \theta_2   [h_{t_1+t+1}(x) -h_{t_1+t}(x)]e^{ \lambda   t }  t^{\frac{1}{\alpha}  -H_t }  \Big|^{\alpha }   dx \nonumber \\
&=&\int_{t_1}^{t_1+1} \big|\theta_2 e^{-\lambda   (t_1-x)} \big( 1 -e^{-\lambda  }   \big)  \big|^{\alpha}   dx  .\label{fmgdg}
\end{eqnarray}
Then (\ref{ineq33}),  (\ref{fgh2h}) and (\ref{fmgdg}) implies that for large $t$,
\begin{eqnarray} \label{inep25}
- I(\theta_1, \theta_2, t_1, t_1+t) \geq -I_4( t ) =|I_4( t )| \geq \frac12 C_3 \,  e^{ -\lambda \alpha t }  t^{\alpha H_t -1},
\end{eqnarray}
where $C_3=\int_{t_1}^{t_1+1} \big|\theta_2 e^{-\lambda   (t_1-x)} \big( 1 -e^{-\lambda  }   \big)  \big|^{\alpha}   dx $ does not depend  on $t.$ Combining (\ref{fh963h}) and (\ref{inep25}) together, we have
 \begin{eqnarray} \label{fv26s}
 I(\theta_1, \theta_2, t_1, t_1+t) \asymp   e^{ -\lambda \alpha t }  t^{ \alpha H_t -1}
\end{eqnarray}
for   $\theta_1 \theta_2 \neq 0.$
It is easy to see that
\begin{eqnarray}
 &&  K(\theta_1, \theta_2, t_1, t_1+t) \nonumber  \\
 &&= \exp\bigg\{ - \int_{-\infty}^{1}\Big|   \theta_1  [e^{-\lambda (1-u)_+}(1-u)_+^{H_{t_1+1} -\frac1{\alpha}}  - e^{-\lambda (-u)_+}(-u)_+^{H_{t_1} -\frac1{\alpha}} ] \Big|^{\alpha}   du\bigg\} \nonumber  \\
 && \ \ \ \times \exp\bigg\{   -\int_{-\infty}^{1}\Big|   \theta_2  [e^{-\lambda ( 1-u)_+}(1-u)_+^{H_{t_1+t+1} -\frac1{\alpha}}  - e^{-\lambda (-u)_+}(-u)_+^{H_{t_1+t} -\frac1{\alpha}} ]  \Big|^{\alpha}   du \bigg\}  \nonumber \\
 && \geq  \exp\bigg\{ -  2 \Big( |\theta_1|^\alpha + |\theta_2|^\alpha \Big)\int_{-\infty}^{1} T(u)    du \bigg\}, \nonumber
 \end{eqnarray}
 where $$T(u):=   e^{-\lambda \alpha(1-u)_+} \Big((1-u)_+^{  \alpha a  -1} + (1-u)_+^{  \alpha b -1} \Big) + e^{-\lambda \alpha (-u)_+}  \Big(( -u)_+^{ \alpha a -1} + ( -u)_+^{\alpha b -1} \Big)  $$
 is integrable on  $ (-\infty, 1]$ with respect to $u$,
and that $ |K(\theta_1, \theta_2, t_1, t_1+t)| \leq 1$.
\, Since $I(\theta_1, \theta_2, t_1, t_1+t)   \rightarrow 0 $ as $t\rightarrow \infty, $
it follows   that $R_{t_1}(t) \sim -K(\theta_1, \theta_2,  t_1, t_1+t)I(\theta_1, \theta_2, t_1, t_1+t);  $
hence (\ref{ineq31}) follows by (\ref{fv26s}).   \qed

\vspace{0.2cm}

\begin{proposition}\label{th42}
Let   $H_t \in [a, b] $  be a continuous function  on $\mathbf{R}.$
Let $Y_{H_t, \alpha , \lambda}(t)$ be a tempered multifractional stable noise (\ref{ineq30}). Recall $R_{t_1}(t)$ defined by (\ref{rtdef})
 with $Y(t)=Y_{H_t, \alpha, \lambda}(t)$. Assume $ \lambda >0,$   $1< \alpha \leq 2$
and $1/\alpha<  H_t  .$
Then
\begin{equation}\label{ineq40}
R_{t_1}(t) \asymp e^{-\lambda  t }t^{H_t-\frac1 \alpha }
\end{equation}
  for  $\theta_1 \theta_2 \neq 0.$
\end{proposition}
\emph{Proof.} Recall $ I_3(t )$  and $ I_4( t) $ defined by (\ref{ineq33}).
By an argument similar to  (\ref{fghlmm5}),  we have for large $t,$
\begin{eqnarray}\label{fsg9n53}
 I_3( t)   I_4( t) \geq 0 .
 \end{eqnarray}
First, we give an estimation for $ I_3(t).$ Using the inequality (\ref{ghllfgm}), we obtain
\begin{eqnarray*}
&& | I_3( t ) |  \leq   2 \int_{-\infty}^{t_1}\alpha \Big|   \theta_2  [h_{t_1+t+1}(x) -h_{t_1+t}(x)]   \Big| \Big| \theta_1 [h_{t_1+1}(x) -h_{t_1}(x)] \Big|^{\alpha-1}   dx \\
&&\ \ \ \   \ \   \ \leq 4| \theta_2| |\theta_1|^{\alpha-1}   \int_{-\infty}^{t_1}\Big| h_{t_1+t+1}(x) -h_{t_1+t}(x)  \Big|   \Big|  h_{t_1+1}(x) -h_{t_1}(x)  \Big|^{\alpha-1}  dx.
\end{eqnarray*}
By an argument similar to (\ref{fine524}), it is easy to see that
for large $t$ and $x\leq t_1$,
\begin{eqnarray}
 \Big|   [h_{t_1+t+1}(x) -h_{t_1+t}(x)]e^{ \lambda   t }  t^{\frac{1}{\alpha}  -H_t } \Big|
 \leq  2 e^{-\lambda   (t_1-x)} (1 + e^{-\lambda  }) (2+ t_1-x  )^{H_t-\frac 1{\alpha}} . \label{ineq43}
\end{eqnarray}
Therefore, for  large $t$ and $x\leq t_1$,
\begin{eqnarray}
|I_3( t )|&\leq&  8| \theta_2| |\theta_1|^{\alpha-1}  \nonumber \\
 &&   \times\int_{-\infty}^{t_1} e^{ -\lambda   t }  t^{H_t-\frac1{\alpha}} (1 +e^{-\lambda  }) \big(2 +t_1-x\big) ^{H_t-\frac{1}{\alpha}  } e^{-\lambda  (t_1-x)}  \Big|  h_{t_1+1}(x) -h_{t_1}(x)  \Big|^{\alpha-1}  dx \nonumber \\
 &\leq&  8| \theta_2| |\theta_1|^{\alpha-1}  e^{ -\lambda   t }  t^{H_t-\frac1{\alpha}}\nonumber \\
 &&   \times\int_{-\infty}^{t_1} (1 +e^{-\lambda  }) \Big(2 +t_1-x\Big) ^{b-\frac{1}{\alpha}  } e^{-\lambda  (t_1-x)}  \Big|  h_{t_1+1}(x) -h_{t_1}(x)  \Big|^{\alpha-1}  dx  .  \label{ineq44}
\end{eqnarray}
From (\ref{ineq44}), we obtain
\begin{eqnarray}\label{ineq46}
   |I_3( t  )|  \leq  C_1 | \theta_2| |\theta_1|^{\alpha-1}  e^{ -\lambda   t }  t^{ H_t-\frac1{\alpha}} ,
\end{eqnarray}
where $C_1$ does not depend  on $t.$ Similarly, we have  for large $t,$
\begin{eqnarray}
|I_4( t)|
 &\leq&   C_2 | \theta_2||\theta_1|^{\alpha-1}  e^{-\lambda t}  t^{H_t-\frac1{\alpha}} ,  \label{ineq47}
\end{eqnarray}
where $C_2$ does not depend  on $t.$
Therefore, from   (\ref{ineq46}) and (\ref{ineq47}), for large $t,$
\begin{eqnarray}\label{fs53f}
|I(\theta_1, \theta_2, t_1, t_1+t)| \leq C_3 | \theta_2||\theta_1|^{\alpha-1} e^{-\lambda t}  t^{H_t-\frac1{\alpha}} .
\end{eqnarray}
where $C_3$ does not depend  on $t.$
  By (\ref{ghllfgm})  we have for large $t,$
\begin{eqnarray*}
  | I_4( t)| \geq \frac12 \int_{t_1}^{t_1+1}\alpha \Big|   \theta_2  [h_{t_1+t+1}(x) -h_{t_1+t}(x)]   \Big| \Big| \theta_1  h_{t_1+1}(x)   \Big|^{\alpha -1}   dx .
\end{eqnarray*}
Applying (\ref{fmcwgdg}) with $\alpha=1$ and the dominated convergence theorem yields
\begin{eqnarray*}
& &\liminf_{t\rightarrow \infty} | I_4( t)| e^{\lambda  t }t^{\frac1 \alpha -H_t}   \nonumber\\
 &&\geq  \lim_{t\rightarrow \infty} \frac12\alpha | \theta_2||\theta_1|^{\alpha-1}  \int_{t_1}^{t_1+1} \Big|   [h_{t_1+t+1}(x) -h_{t_1+t}(x)]e^{ \lambda   t }  t^{\frac{1}{\alpha }  -H_t } \Big|  \Big|  h_{t_1+1}(x)    \Big|^{\alpha -1}  dx \nonumber \\
 &&= \frac12\alpha | \theta_2||\theta_1|^{\alpha-1} \int_{t_1}^{t_1+1} e^{-\lambda   (t_1-x)} \big(1-e^{-\lambda  })  \Big|  h_{t_1+1}(x)    \Big|^{\alpha -1}  dx  .
\end{eqnarray*}
Thus
\begin{eqnarray}\label{fgjvvm5}
  |I_4( t) |\succeq  e^{-\lambda  t }t^{H_t-\frac1 \alpha}   .
\end{eqnarray}
for $\theta_1 \theta_2 \neq 0.$
Combining (\ref{fsg9n53}), (\ref{fs53f}) and (\ref{fgjvvm5}) together, we have
 \begin{eqnarray}
e^{-\lambda  t }t^{H_t-\frac1 \alpha} \preceq |I_4( t)|  \preceq  |I(\theta_1, \theta_2, t_1, t_1+t)|  \preceq   e^{-\lambda  t }t^{H_t-\frac1 \alpha}
\end{eqnarray}
for   $\theta_1 \theta_2 \neq 0.$
Since $I(\theta_1, \theta_2, t_1, t_1+t)   \rightarrow 0 $ as $t\rightarrow \infty, $
it follows from (\ref{rtdef}) that $R_{t_1}(t) \sim -K(\theta_1, \theta_2, t_1, t_1+t)I(\theta_1, \theta_2, t_1, t_1+t);  $
hence (\ref{ineq40}) holds.   \qed

\begin{remark}
One says that a symmetric $\alpha$-stable process $X(t)$ exhibits long-range dependence if for any $t_1\in \mathbf{R}$,
\begin{eqnarray}\label{gsgfhl}
\sum_{n=0}^\infty \Big| R_{t_1}(n)\Big|= \infty,
\end{eqnarray}
where $R_{t_1}(t)$ is defined by (\ref{rtdef}).
 It is obvious that LTFmSM and LTmFSM   are not long-range dependent,
but they exhibit semi-long-range dependence,
that is, for $\lambda>0$ sufficiently small, the sum (\ref{gsgfhl}) is large, and it tends
to infinity as $\lambda\rightarrow 0.$ Therefore, LTFmSM and LTmFSM  provide  two useful alternative models
for data that exhibit strong dependence.
\end{remark}

\section{Scaling property and tail probabilities}\label{sec4}
The following result shows that  LTmFSM (\ref{defin01}) has a nice   scaling property, involving both the time scale and the tempering.
Denote by $\stackrel{fdd}{=} $  equality in the sense of finite dimensional distributions.
\begin{proposition} For any scale factor $c>0$, it holds
\begin{eqnarray}\label{ineq1}
\Big\{ X_{H_{ct},\alpha,\lambda }(c\, t) \Big\}_{t \in \mathbf{R}} \stackrel{fdd}{=} \Big \{ c^{H_{c\,t} } X_{H_{ct},\alpha,c\lambda }(t) \Big\}_{t \in \mathbf{R}}.
\end{eqnarray}
\end{proposition}
\emph{Proof.} It is easy to see that  $$G_{H_{ct},\alpha,\lambda}(c\,t, c\, x) =c^{H_{ct}-\frac1{\alpha}} G_{H_{ct},\alpha,c\lambda }(t, x).$$
 Notice that $dZ_{\alpha}(c\,x)$ has control measure $c^{\frac1 \alpha}dx.$ Given $t_1 < t_2 < ...< t_n,$ a change  of variable
$x=c\,x'$ then yields
\begin{eqnarray*}
 (  X_{H_{c\,t_i},\alpha ,\lambda }(c\,t_i) : i=1,...,n)&=& \Big( \int_{-\infty}^{\infty} G_{H_{c\,t_i},\alpha,\lambda }( c\, t_i, x) dZ_{\alpha}(x) : i=1,...,n \Big) \\
 &=&  \Big( \int_{-\infty}^{\infty} G_{H_{ct_i},\alpha,\lambda }(c\, t_i,   c x')dZ_{\alpha}(c\,x') : i=1,...,n \Big) \\
  &\stackrel{d}{=}&  \Big( \int_{-\infty}^{\infty}c^{H_{ct_i}-\frac1{\alpha}} G_{H_{ct_i},\alpha,c\lambda }( t_i, x')\, c^{\frac1 \alpha}dZ_{\alpha}(x') : i=1,...,n \Big) \\
   &=&  \Big(c^{H_{ct_i}} \int_{-\infty}^{\infty} G_{H_{ct_i},\alpha,c\lambda }(t_i, x') dZ_{\alpha}(x') : i=1,...,n \Big) \\
    &=&  \Big( c^{H_{ct_i}}X_{H_{ct_i},\alpha,c\lambda }(t_i)   : i=1,...,n \Big) ,
\end{eqnarray*}
where $\stackrel{d}{=} $ indicates equality in distribution.
So that (\ref{ineq1}) holds.
   \hfill\qed

We say that a stochastic process $X(t),\, t \in  I,$ is  \emph{stochastic H\"{o}lder continuous} of exponent $\beta \in (0, \infty)$ if it holds
\[
\limsup_{ t,v \in I,\  |t-v| \rightarrow   0 }   \mathbf{P}\big(|X(t)-X(v)|\geq C |t-v|^\beta \big)=0
\]
for a positive constant $C.$  It is obvious  that
if $X(u)$ is  stochastic  H\"{o}lder continuous of exponent $\beta_1 ,$ then $X(u)$ is
stochastic  H\"{o}lder continuous  of exponent $\beta_2 \in (0, \beta_1].$

The following proposition shows that   LTFmSM is stochastic H\"{o}lder continuous. Denote
  $a\wedge b =\min\{a, b\}.$

\begin{proposition}\label{pr2.0}
There is a number $C,$ depending only on $a, b, \lambda$ and $H$, such that for all $t, v \in \mathbf{R}$ and any $y > 0,$
\begin{eqnarray}\label{ite1}
\mathbf{P}\Big( \Big| X_{H ,\alpha(x),\lambda }( t) - X_{H ,\alpha(x),\lambda }( v)\Big|\geq y \Big) \leq\frac{ C }{ y^a \wedge y^b }  \Big(  |t  -v|^{Ha} +  |t -v|^{ Hb}  \Big )
\end{eqnarray}
In particular, (\ref{ite1}) implies that for any $\beta \in (0, Ha/b)$ and all $t, v$ satisfying $|t-v|\leq 1,$
\[
\mathbf{P}\Big( \Big| X_{H ,\alpha(x),\lambda }( t) - X_{H ,\alpha(x),\lambda }( v)\Big|\geq |t-v|^\beta \Big) \leq   C     |t -v|^{Ha-\beta b},
\]
which implies that  $ X_{H ,\alpha(x),\lambda }( t)$  is stochastic H\"{o}lder continuous of exponent $\beta \in (0,  Ha/b )$.
\end{proposition}
\emph{Proof.}
 By Proposition 2.3 of  Falconer and  Liu  \cite{FL12}, it follows  that for any $y>0,$
\begin{eqnarray}
 &&\mathbf{P}\Big( \Big| X_{H ,\alpha(x),\lambda }( t) - X_{H ,\alpha(x),\lambda }( v)\Big|\geq y\Big) \nonumber\\
¡§
&& \ \ \leq   C_1   \int_{-\infty}^{\infty}  \Big|  \frac{G_{H ,\alpha(x),\lambda }(t, x)-G_{H ,\alpha(x),\lambda }(v, x) }{y}  \Big|^{\alpha(x)} dx \nonumber \\
&& \ \ \leq   \frac{ C_1}{y^a \wedge y^b }    \int_{-\infty}^{\infty}  \Big|  G_{H ,\alpha(x),\lambda }(t, x)-G_{H ,\alpha(x),\lambda }(v, x)   \Big|^{\alpha(x)} dx , \label{ineaf}
\end{eqnarray}
where $G_{H ,\alpha(x),\lambda }(t, x)$ is defined by (\ref{funcgs}).
Without loss of generality, we  assume that $t\geq v.$ Then
\begin{eqnarray}\label{ineq61}
 \int_{-\infty}^{\infty}  \Big|  G_{H ,\alpha(x),\lambda }(t, x)-G_{H ,\alpha(x),\lambda }(v, x)   \Big|^{\alpha(x)} dx = I_1 + I_2,
\end{eqnarray}
where
\begin{eqnarray}
 I_1 &=&   \int_{-\infty}^{v}  \Big|  G_{H ,\alpha(x),\lambda }(t, x)-G_{H ,\alpha(x),\lambda }(v, x)   \Big|^{\alpha(x)} dx, \nonumber \\
  I_2  &=&  \int_{v}^{t }   e^{-\lambda \alpha(x) (t-x) }(t-x) ^{H \alpha(x)-1  }  dx . \nonumber
\end{eqnarray}
Using the inequality $|x+ y|^\alpha \leq 2^\alpha (|x|^\alpha + |y|^\alpha)$ for all $x, y \in \mathbf{R}$ and any $\alpha>0,$
we have $$I_1 \leq 4 (I_{11} + I_{12}), $$  where
\begin{eqnarray}
 I_{11} &=&   \int_{-\infty}^{v}  \Big| e^{-\lambda (t-x)}(t-x)^{H -\frac1{\alpha(x)}}-e^{-\lambda (t-x)}(v-x)^{H -\frac1{\alpha(x)}}  \Big|^{\alpha(x)} dx, \nonumber \\
 I_{12} &=&   \int_{-\infty}^{v}  \Big| e^{-\lambda (t-x)}(v-x)^{H -\frac1{\alpha(x)}}-e^{-\lambda (v-x)}(v-x)^{H -\frac1{\alpha(x)}}  \Big|^{\alpha(x)} dx. \nonumber
\end{eqnarray}
Let $h=t-v.$ We deduce the following  estimation of $ I_{11}:$
\begin{eqnarray}
 I_{11} &\leq&    \int_{-\infty}^{v}   \Big|(t -x)^{H -\frac1{\alpha(x)}}-(v-x)^{H -\frac1{\alpha(x)}}  \Big|^{\alpha(x)} dx  \nonumber \\
 &\leq&    \int_{-\infty}^{v}   \Big|(h +v -x)^{H -\frac1{\alpha(x)}}-(v-x)^{H -\frac1{\alpha(x)}}  \Big|^{\alpha(x)} dx  \nonumber \\
 &=&    \int_{-\infty}^{v}   \Big|\Big(1+ \frac{v -x}{h}\Big)^{H -\frac1{\alpha(x)}}-\Big(\frac{v -x}{h}\Big)^{H -\frac1{\alpha(x)}}  \Big|^{\alpha(x)} h^{H \alpha(x)-1}dx  \nonumber \\
  &=&    \int_{0}^{\infty }   \Big|\big(1+ u\big)^{H -\frac1{\alpha(v-hu)}}- u^{H -\frac1{\alpha(v-hu)}}  \Big|^{\alpha(v-hu)} h^{H \alpha(v-hu) }du  \nonumber \\
  &\leq&  \int_{0}^{\infty }   \Big|\big(1+ u\big)^{H -\frac1{\alpha(v-hu)}}- u^{H -\frac1{\alpha(v-hu)}}  \Big|^{\alpha(v-hu)}  du \  \Big(  h^{Ha} +  h^{Hb}  \Big ) \nonumber \\
   &\leq&  C_{11} \Big(  |t -v|^{Ha} +  |t -v|^{Hb}  \Big )  .  \nonumber
\end{eqnarray}
Next, we estimate $ I_{12}.$  Notice that $|e^{-x}-e^{-y}| \leq |x-y|$ for $x, y >0.$ Substitute $u=v-x$
to see that for $\lambda>0,$
\begin{eqnarray}
 I_{12} &=&    \int_{-\infty}^{v} \frac{1}{(\lambda \alpha(x))^{H \alpha(x) -  1} }  \big(\lambda \alpha(x)(v-x)\big)^{H \alpha(x) -  1}  e^{-\lambda \alpha(x)(v-x)}  \Big| e^{-\lambda (t-v)} - 1 \Big|^{\alpha(x)} dx  \nonumber \\
 &\leq&   C_{12} \int_{-\infty}^{v} \Big(\lambda \alpha(x)(v-x)\Big)^{H \alpha(x) -  1}  e^{-\lambda \alpha(x)(v-x)}   \min\Big\{ (t-v)^{\alpha(x)},   1  \Big \}dx  \nonumber \\
 &\leq&  C_{12}\min\Big\{ |t -v|^{a},   1  \Big \} \int_{0}^{\infty}(\lambda \alpha(v-u)u)^{H \alpha(v-u) -  1}  e^{-\lambda \alpha(v-u)u}   du  \nonumber \\
  &\leq&  C_{12}\min\Big\{ |t -v|^{a},   1  \Big \} \int_{0}^{\infty}  \max_{\alpha \in [a, b]}\Big\{ (\lambda \alpha u)^{H \alpha -  1}  e^{-\lambda \alpha u} \Big\} du  \nonumber \\
   &\leq&  C_{13}\min\Big\{ |t-v|^{a},   1  \Big \}  . \label{fds01}
\end{eqnarray}
It is obvious that if  $\lambda=0,$ then $I_{12}=0,$ and thus (\ref{fds01}) holds obviously for all $\lambda\geq 0$.
By simple calculations, we get
\begin{eqnarray*}
I_2  &\leq&   \int_{v}^{t }   (t-x) ^{H \alpha(x)-1  }  dx  \nonumber \\
&\leq& \left\{ \begin{array}{ll}
\int_{v}^{t }   (t-x) ^{H a-1  }  dx & \textrm{if $t -v \leq 1$}\\
\\
 \int_{v}^{t-1 }   (t-x) ^{H b -1  }  dx + \int_{t-1}^{t }   (t-x) ^{H a-1  }  dx\ \ \ \ & \textrm{if $t -v > 1$}
\end{array} \right.  \\
 &\leq&  C_4 \Big(  |t -v|^{Ha} +  |t -v|^{Hb}  \Big ) .
\end{eqnarray*}
Returning to (\ref{ineq61}), we obtain
\begin{eqnarray}
 && \int_{-\infty}^{\infty}  \Big|  G_{H ,\alpha(x),\lambda }(t, x)-G_{H ,\alpha(x),\lambda }(v, x)   \Big|^{\alpha(x)} dx  \nonumber \\
& &\  \leq  C_5 \Big(  |t -v|^{Ha} +  |t -v|^{Hb} + \min\Big\{ |t-v|^{a},   1  \Big \}  \Big )\nonumber \\
& &\  \leq  C_6 \Big(  |t -v|^{Ha} +  |t -v|^{Hb}  \Big ). \label{sjkns}
\end{eqnarray}
Hence, for $y>0,$
\[
 \mathbf{P}\Big( \Big| X_{H ,\alpha(x),\lambda }( t) - X_{H ,\alpha(x),\lambda }( v)\Big|\geq y \Big)
  \leq   \frac{ C_7}{y^a \wedge y^b }  \Big(  |t -v|^{Ha} +  |t -v|^{Hb}  \Big ).
\]
This completes the proof of Proposition \ref{pr2.0}. \hfill\qed

The following proposition shows that   LTmFSM is also stochastic H\"{o}lder continuous.

\begin{proposition}\label{1pr2.0} Let $\lambda>0.$
There is a number $C$ depending only on $a, b$, $\alpha$ and $\lambda,$ such that for all $z > 0,$
\begin{eqnarray}\label{ite}
\mathbf{P}\Big( \Big|  X_{H_t ,\alpha,\lambda }( t) -X_{H_s,\alpha,\lambda }( s)\Big|\geq z \Big) \leq \frac{ C}{ z^\alpha }  \Big(  |t -s|^{\alpha H_t} + \big|{H_t  }- {H_s } \big|^{\alpha} \Big)
\end{eqnarray}
for all $t, s \in \mathbf{R}$ satisfying $t\geq  s$.
In particular, if $H_t$ is $\gamma-$H\"{o}lder continuous, that is
$$\big|{H_t  }- {H_s } \big|  \leq C  |t-s|^\gamma  \ \ \ \  \textit{for}\  t, s \ \textit{satisfying} \  |t-s| \leq 1,$$
 then
(\ref{ite}) implies that for any $\beta \in (0,  \min\{a,  \gamma\})$ and all $t, s \in \mathbf{R}$ satisfying $|t-s| \leq  1,$
\[
\mathbf{P}\Big( \Big| X_{H_t ,\alpha ,\lambda }( t) - X_{H_s ,\alpha ,\lambda }( s)\Big|\geq |t-s|^\beta \Big) \leq   C \,  \Big(  |t -s|^{\alpha (a- \beta) } + |t -s|^{ \alpha (\gamma -   \beta) } \Big),
\]
which implies that  $ X_{H_t ,\alpha,\lambda }( t)$  is stochastic H\"{o}lder continuous of exponent $\beta \in (0,  \min\{a,  \gamma\})$.
\end{proposition}
\emph{Proof.} By Proposition 1.2.15 of Samorodnitsky and  Taqqu \cite{ST94},  it follows  that for $z>0,$
\begin{eqnarray}
 &&  \mathbf{P}\Big( \Big| X_{H_t ,\alpha,\lambda }( t) - X_{H_s ,\alpha,\lambda }( s)\Big|\geq z\Big) \nonumber\\
 &&\ \ \ \  \leq  C_0  \frac{ 1}{z^\alpha } \int_{-\infty}^{\infty}  \Big|G_{H_t ,\alpha,\lambda }(t, x)-G_{H_s ,\alpha,\lambda }(s, x) \Big|^{\alpha} dx  . \label{ineaf}
\end{eqnarray}
Using the inequality  for any $\alpha>0,$
$$|x+ y + z |^\alpha \leq 3^\alpha (|x|^\alpha + |y|^\alpha + |z|^\alpha), \ \ \ \ \ \ \  x, y, z \in \mathbf{R},$$
we have
\begin{eqnarray}\label{ineq6}
 \int_{-\infty}^{\infty}  \Big|  G_{H_t ,\alpha,\lambda }(t, x)-G_{H_s ,\alpha,\lambda }(s, x)   \Big|^{\alpha} dx \leq 3^\alpha ( I_1 + I_2 +I_3 ),
\end{eqnarray}
where
\begin{eqnarray}
 I_1 &=&   \int_{-\infty}^{\infty}  \Big|  e^{-\lambda (t-x)_+}(t-x)_+^{H_t -\frac1{\alpha}}-    e^{-\lambda (s-x)_+}(s-x)_+^{H_t -\frac1{\alpha}}   \Big|^{\alpha} dx, \nonumber \\
  I_2 &=&   \int_{-\infty}^{\infty}  \Big|  e^{-\lambda (s-x)_+}(s-x)_+^{H_t -\frac1{\alpha}}-    e^{-\lambda (s-x)_+}(s-x)_+^{H_s -\frac1{\alpha}}   \Big|^{\alpha} dx, \nonumber \\
  I_3  &=&  \int_{-\infty}^{\infty}   \Big|  e^{-\lambda (-x)_+}(-x)_+^{H_t -\frac1{\alpha}}-    e^{-\lambda (-x)_+}(-x)_+^{H_s -\frac1{\alpha}}   \Big|^{\alpha} dx . \nonumber
\end{eqnarray}
It is easy to see that
\begin{eqnarray}\label{hjls2}
I_1 \leq 2^\alpha (I_{11} + I_{12}),
\end{eqnarray}
 where
\begin{eqnarray}
 I_{11} &=&   \int_{-\infty}^{\infty}  \Big| e^{-\lambda (t-x)_+}(t-x)_+^{H_t -\frac1{\alpha}}-e^{-\lambda (t-x)_+}(s-x)_+^{H_t -\frac1{\alpha}}  \Big|^{\alpha} dx, \nonumber \\
 I_{12} &=&   \int_{-\infty}^{\infty}  \Big| e^{-\lambda (t-x)_+}(s-x)_+^{H_t -\frac1{\alpha}}-e^{-\lambda (s-x)_+}(s-x)_+^{H_t -\frac1{\alpha}}  \Big|^{\alpha} dx. \nonumber
\end{eqnarray}
Let $h=t-s >0.$ Notice that $|\big(1+ u\big)^{H_t -\frac1{\alpha}}- u^{H_t -\frac1{\alpha}} | \leq 2 \beta u^{H_t -\frac1{\alpha}-1} , u\rightarrow \infty.$ Then we deduce the following  estimation of $ I_{11}:$
\begin{eqnarray}
 I_{11} &=&    \int_{-\infty}^{t}e^{-\lambda \alpha(t-x)}  \Big|(t -x)^{H_t -\frac1{\alpha}}-(s-x)_+^{H_t -\frac1{\alpha}}  \Big|^{\alpha} dx  \nonumber \\
 &=&    \int_{-\infty}^{t} e^{-\lambda \alpha(t-x)}  \Big|\Big(1+ \frac{s -x}{h}\Big)^{H_t -\frac1{\alpha}}-\Big(\frac{s -x}{h}\Big)_+^{H_t -\frac1{\alpha}}  \Big|^{\alpha} h^{H_t \alpha-1}dx  \nonumber \\
  &=&    \int_{-1}^{\infty }  e^{-\lambda \alpha h (1+u)}   \Big|\big(1+ u\big)^{H_t -\frac1{\alpha}}- u_+^{H_t-\frac1{\alpha}}  \Big|^{\alpha} h^{H_t \alpha   }du  \nonumber \\
  &\leq&  \int_{-1}^{\infty }  e^{-\lambda \alpha h (1+u)}   \Big|\big(1+ u\big)^{H_t -\frac1{\alpha}}- u_+^{H_t -\frac1{\alpha}}  \Big|^{\alpha}  du \   h^{H_t \alpha   } \nonumber \\
    &\leq& C_{11}  h^{H_t \alpha   }= C_{11}  |t -s|^{H_t \alpha   } .  \label{ssggjf}
\end{eqnarray}
Next, consider the  item  $ I_{12}.$ Substitute $u=s-x$ and then $w=\lambda$
to see that for $\lambda>0,$
\begin{eqnarray}
 I_{12}
   &=& \int_{-\infty}^{s} (s-x)^{\alpha H_t -1}    e^{-\lambda \alpha(s-x)}  \Big| e^{-\lambda (t-s)} - 1 \Big|^{\alpha} dx  \nonumber \\
 &\leq&  \int_{-\infty}^{s}  (s-x)^{\alpha H_t -  1}   e^{-\lambda \alpha(s-x)}   dx  \, \min\Big\{ (t-s)^{\alpha},   1  \Big \} \nonumber \\
 &=&   \int_{0}^{\infty}u^{\alpha H_t -  1}   e^{-\lambda \alpha u}   du  \, \min\Big\{ (t-s)^{\alpha},   1  \Big \} \nonumber \\
   &\leq&  C_{12}  \min\Big\{ |t-s|^{\alpha},   1  \Big \},  \label{tnjs}
\end{eqnarray}
where the second line of the last inequalities follows by  the inequality $|e^{-x}-e^{-y}| \leq |x-y|$ for all $x, y \geq 0.$
It is obvious that if  $\lambda=0,$ then $I_{12}=0.$ Thus (\ref{tnjs}) also holds for $\lambda=0$.
Combining (\ref{hjls2}), (\ref{ssggjf}) and (\ref{tnjs}) together,  we get
\begin{eqnarray}\label{ineqs12}
I_1 \leq C_1 \Big( |t -s|^{\alpha  H_t } + \min\Big\{ |t-s|^{\alpha},   1  \Big \} \Big).
\end{eqnarray}

In the sequel, we give the estimations of $I_2$ and $I_3$.
Without loss of generality, we  assume that $H_t\geq  H_s.$
By some simple calculations, we get
\begin{eqnarray}
  I_2
  &=&   \int_{-\infty}^{s} e^{-\lambda \alpha (s-x) } \Big|(s-x) ^{H_t -\frac1{\alpha}}-    (s-x) ^{H_s -\frac1{\alpha}}   \Big|^{\alpha} dx \nonumber \\
  &=&   \int_{0}^{\infty } e^{-\lambda \alpha u } u^{-1}\Big|u ^{H_t }-    u^{H_s  }   \Big|^{\alpha} du \nonumber \\
  &=& \int_{0}^{\infty} e^{-\lambda \alpha u}  \Big|{H_t  }- {H_s } \Big|^{\alpha}  u^{\alpha H_\theta -1} |\log u | ^ \alpha du \nonumber \\
  &\leq&   C_2 \big|{H_t  }- {H_s } \big|^{\alpha}, \label{ineqsqs10}
\end{eqnarray}
where $H_\theta \in [H_s, H_t].$ Similarly, we have
\begin{eqnarray}
I_3  &\leq&    C_3 \big|{H_t  }- {H_s } \big|^{\alpha}. \label{ineq10}
\end{eqnarray}
Combining the inequalities (\ref{ineq6}), (\ref{ineqs12}), (\ref{ineqsqs10}) and (\ref{ineq10}) together,  we obtain
\begin{eqnarray}
&&  \int_{-\infty}^{\infty}  \Big|  G_{H_t ,\alpha,\lambda }(t, x)-G_{H_s ,\alpha,\lambda }(s, x)   \Big|^{\alpha} dx\nonumber \\
& &\leq  C_4 \Big(  |t -s|^{\alpha H_t}  +  \min\big\{ |t-s|^{\alpha},   1  \big \}  + \big|{H_t  }- {H_s } \big|^{\alpha} \Big )\nonumber \\
& &\leq  C_5 \Big(  |t -s|^{\alpha H_t}   + \big|{H_t  }- {H_s } \big|^{\alpha} \Big ). \label{1sjkns}
\end{eqnarray}
Returning to (\ref{ineaf}), we get  for $z>0,$
\begin{eqnarray*}
 \mathbf{P}\Big( \Big| X_{H_t ,\alpha,\lambda }( t) - X_{H_s ,\alpha,\lambda }( s)\Big|\geq z \Big)
  \leq   \frac{C_6 }{z^\alpha }  \Big(   |t -s|^{\alpha H_t} + \big|{H_t  }- {H_s } \big|^{\alpha}  \Big ).
\end{eqnarray*}
This completes the proof of Proposition \ref{1pr2.0}. \hfill\qed

\section{Absolute  moments}
We estimate the absolute (incremental) moments of the LTFmSM.
\begin{proposition}\label{pr2.11}
If $0<p< a$, then there exists a number  $C_1,$  depending only on $ a, b, \lambda$ and $H,$ such that for all $t, v \in \mathbf{R}$ and $|t-v|\geq 1$,
\[
\mathbf{E}\Big[\big| X_{H ,\alpha(x),\lambda }( t) - X_{H ,\alpha(x),\lambda }( v)\big|^p\Big] \leq C_1\Big(1+ \frac{p}{a-p} \Big)     |t -v|^{Hb} .
\]
\end{proposition}
\emph{Proof.} When $|t-v|\geq 1,$ using Proposition \ref{pr2.0},   we deduce that
\begin{eqnarray}
&& \mathbf{E}\Big[\big| X_{H ,\alpha(x),\lambda }( t) - X_{H ,\alpha(x),\lambda }( v)\big|^p\Big] \nonumber\\
&& =p\int_0^\infty    y^{p-1} \mathbf{P}\Big( \Big| X_{H ,\alpha(x),\lambda }( t) - X_{H ,\alpha(x),\lambda }( v)\Big|\geq y \Big) dy  \nonumber \\
&& \leq p  \Big(\int_0^1  y^{p-1}   dy +  C_1\int_1^\infty    y^{p-1-a}  dy \, \Big)    |t -v|^{Hb}    \nonumber \\
&& \leq C_2  \Big(1+ \frac{p}{a-p} \Big)    |t -v|^{Hb}   . \nonumber
\end{eqnarray}
This completes the proof of Proposition \ref{pr2.11}. \hfill\qed

The next proposition gives an estimate for the absolute (incremental) moment of the LTmFSM.
\begin{proposition}\label{1pr2.11}
If $0<p< \alpha$, then there is a number $C$ depending only on $p, a, b$ and $ \lambda$, such that for all $t, s \in \mathbf{R}$ satisfying $t \geq s$,
\begin{eqnarray}
  \mathbf{E}\Big[\big| X_{H_t ,\alpha ,\lambda }( t) - X_{H_s ,\alpha ,\lambda }( s)\big|^p\Big]
\leq \Big(1+   \frac{2C_1}{p-\alpha} \Big)   \Big(  |t -s|^{p H_t}  + \big|{H_t  }- {H_s } \big|^{p}   \Big ). \nonumber
\end{eqnarray}
\end{proposition}
\emph{Proof.} Using Proposition \ref{1pr2.0}, we have for any $\varepsilon>0,$
\begin{eqnarray}
&& \mathbf{E}\Big[\big| X_{H_t ,\alpha ,\lambda }( t) - X_{H_s ,\alpha ,\lambda }( s)\big|^p\Big] \nonumber \\
&& =p\int_0^\infty    y^{p-1} \mathbf{P}\Big( \Big| X_{H_t ,\alpha ,\lambda }( t) - X_{H_s ,\alpha ,\lambda }( s)\Big|\geq y \Big) dy  \nonumber \\
&& \leq   p \int_0^\varepsilon  y^{p-1}   dy  + C_1  \int_\varepsilon^\infty    y^{p-1-\alpha}  dy   \Big(  |t -s|^{\alpha H_t} + \big|{H_t  }- {H_s } \big|^{\alpha}   \Big )  \nonumber \\
&&= \varepsilon^p+  \frac{C_1}{\alpha -p} \varepsilon^{p-\alpha}   \Big(  |t -s|^{\alpha H_t}  + \big|{H_t  }- {H_s } \big|^{\alpha}   \Big ). \nonumber
\end{eqnarray}
Taking   $ \varepsilon= \max\{ |t -s|^{  H_t} ,  \big| H_t  -  H_s   \big| \} ,$ we get
\begin{eqnarray}
  \mathbf{E}\Big[\big| X_{H_t ,\alpha ,\lambda }( t) - X_{H_s ,\alpha ,\lambda }( s)\big|^p\Big]
\leq \Big(1+   \frac{2C_1}{\alpha-p} \Big)   \Big(  |t -s|^{p H_t}  + \big|{H_t  }- {H_s } \big|^{p}   \Big ), \nonumber
\end{eqnarray}
which gives the desired inequality.  \hfill\qed

For LFmSM, Le Gu\'{e}vel and L\'{e}vy V\'{e}hel \cite{GLL13a} have investigated
 the asymptotic behaviour of $ \mathbf{E}  \big[ |X ( t+r ) - X ( t )|^\eta \big], r\rightarrow0,$ for some positive constant  $\eta>0$.
The following proposition gives a  result similar  to the one of Le Gu\'{e}vel and L\'{e}vy V\'{e}hel for LTFmSM.
\begin{proposition}\label{fsdf} For  each $t \in \mathbf{R} $  satisfying $ H \alpha(t)\neq 1   $ and  all $\gamma \in (0, a),$    it holds
\begin{eqnarray*}
 \lim_{ r\rightarrow 0+}  \frac{ \mathbf{E}  \big[ |X_{H ,\alpha(x),\lambda }( t+r ) - X_{H ,\alpha(x),\lambda }( t )|^\gamma \big]}{r^{\gamma H} }   =F(\gamma, t),
\end{eqnarray*}
where
\begin{eqnarray*}
F(\gamma, t) = \bigg(\int_{-\infty}^{\infty}\Big[  (1-x)_+^{H -\frac1{\alpha(t)}} -  (-x)_+^{H-\frac1{\alpha(t)}}\Big]^{ \alpha(t)}  dx \bigg) ^{\gamma/\alpha(t)} \frac{2^{\gamma-1} \Gamma\Big(1-\frac{\gamma}{ \alpha(t)}\Big)}{\gamma \int_0^{\infty} u^{-\gamma-1} \sin^2(u)du }
\end{eqnarray*}
and $\Gamma(t)=\int_0^\infty x^{t-1}e^{-x} dx$ is the gamma function.
\end{proposition}
 \emph{Proof.}   Notice that  for all $\gamma \in (0, a)$ and all $  u\in [0,  \, 1),$
\begin{eqnarray}
 &&\mathbf{E} \Big[ \Big|\frac{ X_{H ,\alpha(x),\lambda }( t+r ) - X_{H ,\alpha(x),\lambda }( t )}{r^{H} } \Big|^\gamma \Big] \nonumber\\
  &&\ \  =\ \gamma  \int_{0}^\infty
 z^{\gamma-1}  \ \mathbf{P} \Big( \Big|\frac{ X_{H ,\alpha(x),\lambda }( t+r ) - X_{H ,\alpha(x),\lambda }( t )}{r^{H} }
  \Big| \geq z \Big) dz . \nonumber
\end{eqnarray}
Notice that $X_{H ,\alpha(x),\lambda }( t )$  is $H-$localisable  to $X$ defined by (\ref{ffgls}) (cf. Proposition \ref{pr2.1}  whose proof does not involve  Proposition \ref{fsdf}). Thus
\[
\mathbf{P} \bigg( \Big|\frac{ X_{H ,\alpha(x),\lambda }( t+r ) - X_{H ,\alpha(x),\lambda }( t )}{r^{H} } \Big|  \geq z \bigg) \rightarrow \mathbf{P} \bigg( \big| X (1)\big|  \geq z \bigg), \ \ \ \ \   r\rightarrow0.
\]
By Proposition \ref{pr2.0}, for $z$ large enough,
\begin{eqnarray}
\mathbf{P} \Big( \Big|\frac{ X_{H ,\alpha(x),\lambda }( t+r ) - X_{H ,\alpha(x),\lambda }( t )}{r^{H} }
  \Big| \geq z \Big)   \leq   C \frac{1} {z^a   }   . \nonumber
\end{eqnarray}
Hence, by the Lebesgue dominated convergence theorem,  we have
\begin{eqnarray*}
&&\lim_{r \rightarrow 0+}\mathbf{E} \Big[ \Big|\frac{ X_{H ,\alpha(x),\lambda }( t+r ) - X_{H ,\alpha(x),\lambda }( t )}{r^{H} } \Big|^\gamma \Big]  \\
&& \ =  \gamma  \int_{0}^\infty z^{\gamma-1}  \ \mathbf{P} \bigg( \big| X (1)\big|  \geq z \bigg) dz \\
&& \ =  \mathbf{E}  \big[ \big| X (1)\big|^\gamma \big] \\
& &\ = \bigg(\int_{-\infty}^{\infty}\Big[  (1-x)_+^{H -\frac1{\alpha(t)}} -  (-x)_+^{H-\frac1{\alpha(t)}}\Big]^{ \alpha(t)}  dx \bigg) ^{\gamma/\alpha(t)} \frac{2^{\gamma-1} \Gamma\Big(1-\frac{\gamma}{ \alpha(t)}\Big)}{\gamma \int_0^{\infty} u^{-\gamma-1} \sin^2(u)du },
\end{eqnarray*}
where $X$ is $\alpha(t)$-stable.   We refer to Property 1.2.17 of  Samorodnitsky and Taqqu \cite{ST94} for the last line of the last equality.
 \qed

\section{Sample path properties}
 When $Ha> 1$ with $a>1$, the following proposition  implies that every LTFmSM process  has an a.s.\ H\"{o}lder continuous version.
\begin{proposition}\label{pr2.12}
If $Ha >1$ with $a>1$, then  for any $0< \beta < H -1/a,$ $X_{H ,\alpha(x),\lambda }( t)$   has a continuous version such that
its paths  are almost surely  $\beta-$H\"{o}lder continuous on each bounded interval.
\end{proposition}
\emph{Proof.}   Recall that
$$X_{H  ,\alpha(x),\lambda }( t) = \int_{-\infty}^{\infty}G_{H ,\alpha(x),\lambda }(t, x)  dM_{\alpha}(x).$$
By (\ref{sjkns}), we have for $|t-v| \leq 1,$
\begin{eqnarray}
 \int_{-\infty}^{\infty}  \Big|  G_{H ,\alpha(x),\lambda }(t, x)-G_{H ,\alpha(x),\lambda }(v, x)   \Big|^{\alpha(x)} dx
\leq  C_1 |t-v|^{Ha}.
\end{eqnarray}
By Proposition 3.1 of Falconer and  Liu \cite{FL12}, $X_{H ,\alpha(x),\lambda }( t)$   has a continuous version such that
its paths  are almost surely  $\beta-$H\"{o}lder continuous on each bounded interval, where $0< \beta < (Ha-1)/a.$
  \hfill\qed

Recall that a stochastic process $X(t), t \in T,$ on a probability space $(\Omega, \mathcal{F}, \mathbf{P})$
is called separable if there is a countable set $T^*\subset T$ and an even $\Omega_0 \in \mathcal{F}$ with
$\mathbf{P}(\Omega_0)=0$, such that for any closed set $F \subset \mathbf{R} $ we have
$$\{ \omega: X(t) \in F, \ \ \forall t \in T^*\} \setminus \{\omega: X(t) \in F, \ \ \forall t \in T\} \subset \Omega_0.$$
See Chapter 9 of Samorodnitsky  and Taqqu \cite{ST94} for more details.

When $H_t \alpha< 1$ and $\lambda>0,$ the following proposition  shows that
every separable version of   LTmFSM process has unbounded paths.
\begin{proposition}\label{pr2.13}
If $H_t \alpha< 1$ and $\lambda>0$, then for any separable version of the LTmFSM process, we have for any   interval $(c,\, d),$
\[
\mathbf{P}\Big( \Big\{ \omega: \  \sup_{t \in (c,\, d)}\big| X_{H_t ,\alpha ,\lambda }( t) \big|=\infty \Big\} \Big) =1.
\]
\end{proposition}
\emph{Proof.} We may assume that $(c,\, d)$ is bound.
Consider the countable set $T^*:=\mathbf{Q}\cap[c,\, d],$
where $\mathbf{Q}$ denotes the set of rational numbers. Since $T^*$ is dense in $[c,\, d]$, there exists a sequence of numbers
$\{t_n\}_{n\in \mathbf{N}} \in T^*$, such that for any $x \in [c,\, d],$ $t_n \rightarrow x$ as $n \rightarrow \infty.$
Therefore, it holds
$$f^*(T^*; x):=\sup_{t \in T^*}\Big| G_{H_t ,\alpha ,\lambda }(t, x)\Big| \geq \sup_{t_n \in T^*}\Big| G_{H_t ,\alpha ,\lambda }(t_n, x)\Big|=: f_n^*(T^*; x)=\infty, \ \  n\rightarrow \infty.  $$
Thus  $\int_c^d f^*(T^*; x)dx=\infty,$ and this contradicts Condition (10.2.14)
of Theorem 10.2.3 in Samorodnitsky  and Taqqu \cite{ST94}. Therefore, the stochastic process $\{X_{H_t,\alpha ,\lambda}\}$ does not
have a version with bounded paths on the interval $(c,\, d),$ and this completes the proof.
    \hfill\qed

For LTmFSM process with $H_t\alpha> 1,$  we have the following proposition.
\begin{proposition}\label{1pr2.12}
Assume that $H_t$ is $\gamma-$H\"{o}lder continuous, $\gamma > 1/\alpha,$ that is
\begin{eqnarray}
\big|{H_t  }- {H_s } \big|  \leq C  |t-s|^\gamma    \label{sffg1s}
\end{eqnarray}
for $t, s \in \mathbf{R}$ satisfying $|t-s|\leq 1.$
If $\alpha\min\{a, \gamma \}   >1$, then for any $0< \beta < \min\{a, \gamma \} -1/\alpha,$ $X_{H_t ,\alpha,\lambda }( t)$   has a continuous version, such that
its paths  are almost surely  $\beta-$H\"{o}lder continuous on each compact set.
\end{proposition}
\emph{Proof.} By Proposition \ref{1pr2.11} and (\ref{sffg1s}), we have for  any $0<p<\alpha$  and all $t, s$ satisfying $|t-s|\leq 1,$
\begin{eqnarray}
\mathbf{E}\Big[\big| X_{H_t ,\alpha ,\lambda }( t) - X_{H_s ,\alpha ,\lambda }( s)\big|^p\Big] &\leq& C_1  \Big(  |t -s|^{p a} +   \big|{H_t  }- {H_s } \big|^{p}   \Big ) \nonumber  \\
&\leq& C_2  \Big(  |t -s|^{p a} +   \big|t-s \big|^{p \gamma}   \Big )   .\nonumber
\end{eqnarray}
The Kolmogorov continuity theorem   implies that $X_{H_t ,\alpha ,\lambda }( t)$ has a continuous version, such that
its paths  are almost surely  $\beta-$H\"{o}lder continuous  on each compact set, $0< \beta < (   p\min\{a, \gamma \} -1)/p.$ Let $p\rightarrow \alpha.$
We completes the proof of Proposition \ref{1pr2.12}.   \hfill\qed

\begin{remark} For LmFSM,
 Ayache and Hamonier \cite{AH14} have obtained the uniform pointwise H\"{o}lder exponent of $X_{H_t ,\alpha,\lambda }( t)$. By Theorem 8.1 of Ayache and Hamonier \cite{AH14}, it is easy to see that when $a\geq \gamma,$ the $\beta$ in Proposition 13 cannot exceed $\gamma   -1/\alpha$.
\end{remark}

Denote by $$\mathcal{\widetilde{H}}_t(\omega)=\sup\bigg\{ \gamma: \  \lim_{r\rightarrow 0}\frac{| X_{H ,\alpha(x),\lambda }( t+r, \omega) - X_{H ,\alpha(x),\lambda }(t, \omega)  |}{|r|^\gamma} =0 \bigg\}$$
the pontwise  H\"{o}lder exponent of the LTFmSM $X_{H,\alpha(x),\lambda}(\cdot)$ at $t.$
\begin{proposition}\label{pr22.12}
If $Ha >1$  with $a>1$, then  $\mathcal{\mathcal{\widetilde{H}}}_t(\omega)  \geq  H -1/a$ almost surely.
\end{proposition}
\emph{Proof.} It follows by Proposition \ref{pr2.12}. \hfill\qed
Let $$\mathcal{\widehat{H}}_t(\omega)=\sup\bigg\{ \gamma: \  \lim_{r\rightarrow 0}\frac{| X_{H_t,\alpha,\lambda}( t+r, \omega) - X_{H_t,\alpha,\lambda}(t, \omega)  |}{|r|^\gamma} =0 \bigg\}$$
be the pontwise H\"{o}lder exponent of the LTmFSM  $X_{H_t,\alpha,\lambda}( \cdot)$ at $t.$
\begin{proposition}\label{pr22212}
Assume that $H_t$ is $\gamma-$H\"{o}lder continuous, $\gamma > 1/\alpha.$
If $\alpha\min\{H_{t_0}, \gamma \}   >1$ for some $t_0  \in \mathbf{R}$, then $\mathcal{\widehat{H}}_{t_0}(\omega)    \geq  \min\{H_{t_0}, \gamma \} -1/\alpha$ almost surely.
\end{proposition}
\emph{Proof.} Since $H_t$ is continuous,  we have for any $\varepsilon >0,$ there exists a $\delta>0$  such that for all $s \in [t_0-\delta,\ t_0+\delta],$
it holds $H_s \in [H_{t_0}-\varepsilon,\  H_{t_0}+\varepsilon].$
If $\alpha\min\{H_{t_0}, \gamma \}   >1$, by an argument similar to the proof of Proposition \ref{1pr2.12},  then for any $0< \beta < \min\{H_{t_0}-\varepsilon, \gamma \} -1/\alpha,$ $X_{H_s ,\alpha,\lambda }( s)$   has a continuous version, such that
its paths  are almost surely  $\beta-$H\"{o}lder continuous on $s \in [t_0-\delta,\ t_0+\delta]$. Thus if $\alpha\min\{H_{t_0}, \gamma \}   >1$, then $\mathcal{\widehat{H}}_{t_0}(\omega)    \geq  \min\{H_{t_0}-\varepsilon, \gamma \} -1/\alpha$ almost surely.
  The claim follows by the fact that $\varepsilon$ can be arbitrary small.
    \hfill\qed

\section{H\"{o}lder continuity of quasi norm}
Denote by
$$\Big|\Big|X_{H,\alpha(x),\lambda}(t)\Big|\Big|_\alpha :=  \Bigg\{ y>0 :\ \ \int_{-\infty}^{\infty}\Big|\frac{G_{H,\alpha(x),\lambda}(t, x)}{ y}   \Big|^{\alpha(x)} dx =1  \Bigg\}$$
for $t\in \mathbf{R}.$ Then  $\big|\big| \cdot\big|\big|_\alpha$ is a quasi norm.
In particular, if $\alpha(x)\equiv p\geq 1$ for a constant $p$,   then $  ||X_{H,p,\lambda}(t) | |_p $
is the $L^p(\mathbf{R}) $ norm of $G_{H,p ,\lambda}(t, x).$ Moreover, when $\alpha(x)\equiv \alpha$ for a constant $\alpha \in (0, 2],$
then it holds
$$\Big|\Big|X_{H,\alpha,\lambda}(t)\Big|\Big|_\alpha =\Big(- \log  \mathbf{E}\big[e^{i  X_{H,\alpha,\lambda}(t)  } \big] \Big )^{1/\alpha}=\Bigg(\int_{-\infty}^{\infty}\Big|G_{H,\alpha ,\lambda}(t, x)     \Big|^{\alpha } dx \Bigg)^{1/\alpha}, $$
see Meerschaert and Sabzikar \cite{MS16}.

The next proposition implies  that the quasi norm of LTFmSM process is H\"{o}lder continuous in time $t.$
\begin{proposition}\label{hlmsf}  There are two positive   numbers $c$ and $ C,$ depending only on $ a, b, \lambda$ and $H$, such that
\[
 c \,  |t -v|^{Hb/a}   \leq \Big|\Big| X_{H ,\alpha(x),\lambda }( t) - X_{H ,\alpha(x),\lambda }( v)\Big|\Big|_\alpha  \leq C \, |t -v|^{Ha/b}
\]
for all $t, v \in \mathbf{R}$ satisfying $|t-v|\leq 1$.
\end{proposition}
\emph{Proof.} Denote by $\rho=\big|\big| X_{H ,\alpha(x),\lambda }( t) - X_{H ,\alpha(x),\lambda }( v)\big|\big|_\alpha.$
 Assume that $t> v,$ and write
\begin{eqnarray}
&&\int_{-\infty}^{\infty}  \Big|  G_{H ,\alpha(x),\lambda }(t, x)-G_{H ,\alpha(x),\lambda }(v, x)   \Big|^{\alpha(x)} dx  \nonumber \\
&&\geq \int_{v}^{t }   e^{-\lambda \alpha(x) (t-x) }(t-x) ^{H \alpha(x)-1  }  dx   \nonumber \\
&&\geq e^{-\lambda b (t-v) }\int_{v}^{t }(t-x) ^{H \alpha(x)-1  }  dx  \nonumber \\
&&\geq e^{-\lambda b (t-v) }\int_{v}^{t }(t-x) ^{Hb-1  }  dx  \nonumber \\
&&\geq e^{-\lambda b  }\frac1{Hb} (t-v) ^{Hb  } \nonumber
\end{eqnarray}
uniformly for all $t, v \in \mathbf{R}$ satisfying $|t-v|\leq 1$.  Therefore, we have
\begin{eqnarray*}
1&=&\int_{-\infty}^{\infty}  \Big|\frac{ G_{H ,\alpha(x),\lambda }(t, x)-G_{H ,\alpha(x),\lambda }(v, x) }{\rho} \Big|^{\alpha(x)} dx \\
 &\geq&\int_{-\infty}^{\infty}  \Big|G_{H ,\alpha(x),\lambda }(t, x)-G_{H ,\alpha(x),\lambda }(v, x) \Big|^{\alpha(x)} dx \, \min\Big\{\frac{1 }{\rho^a}, \frac{1 }{\rho^b} \Big\}   \\
 &\geq& e^{-\lambda b  }\frac1{Hb} (t-v) ^{Hb  } \, \min\Big\{\frac{1 }{\rho^a}, \frac{1 }{\rho^b} \Big\}.
\end{eqnarray*}
The last inequality implies the lower bound of $\rho.$
By (\ref{sjkns}), we have
\begin{eqnarray*}
 \int_{-\infty}^{\infty}  \Big|  G_{H ,\alpha(x),\lambda }(t, x)-G_{H ,\alpha(x),\lambda }(v, x)   \Big|^{\alpha(x)} dx
\leq  C_1 |t-v|^{Ha  }
\end{eqnarray*}
uniformly for all $t, v \in \mathbf{R}$ satisfying $|t-v|\leq 1$.
Then
\begin{eqnarray}
1&=&\int_{-\infty}^{\infty}  \Big|\frac{ G_{H ,\alpha(x),\lambda }(t, x)-G_{H ,\alpha(x),\lambda }(v, x) }{\rho} \Big|^{\alpha(x)} dx \nonumber \\
 &\leq&\int_{-\infty}^{\infty}  \Big|G_{H ,\alpha(x),\lambda }(t, x)-G_{H ,\alpha(x),\lambda }(v, x) \Big|^{\alpha(x)} dx \  \max\Big\{\frac{1 }{\rho^a}, \frac{1 }{\rho^b} \Big\}  \nonumber \\
 &\leq&   C_1 |t-v| ^{Ha  } \  \max\Big\{\frac{1 }{\rho^a}, \frac{1 }{\rho^b} \Big\} , \label{sf02f}
\end{eqnarray}
whenever $|t-v|\leq 1$.
Inequality (\ref{sf02f}) implies the upper bound of $\rho.$   \hfill\qed

 When  $ a=b$  and $ 1/a < H  < 1,$  Proposition \ref{hlmsf}  reduces to  Lemma 4.2 of   Meerschaert and Sabzikar \cite{MS16}. Hence Proposition \ref{hlmsf} can be regarded as  a generalization of this lemma.

The next proposition implies  that the quasi norm of LTmFSM process is H\"{o}lder continuous in time $t.$
\begin{proposition}\label{1hlmsf}
There exist two positive   numbers  $c$ and $ C,$ depending only on $ a, b, \lambda $ and $\alpha$, such that
\begin{eqnarray}\label{f1d7}
c \, |t -s|^{H_t}    \leq \Big|\Big| X_{H_t ,\alpha ,\lambda }( t) - X_{H_s ,\alpha ,\lambda }(s)\Big|\Big|_\alpha  \leq C \, \Big(  |t -s|^{  H_t}   + \big|{H_t  }- {H_s } \big|  \Big )
\end{eqnarray}
for all $t, s$ satisfying $0\leq s \leq t \leq s+1$.
\end{proposition}
\emph{Proof.} From the poof of Proposition \ref{1pr2.0}, we have
\begin{eqnarray}
 \int_{-\infty}^{\infty}  \Big|  G_{H_t ,\alpha ,\lambda }(t, x)-G_{H_t ,\alpha ,\lambda }(v, x)   \Big|^{\alpha } dx
&\leq&  C_1    \Big(  |t -s|^{\alpha H_t}   + \big|{H_t  }- {H_s } \big|^{\alpha} \Big )   \nonumber \\
&\leq&  2C_1   \max \Big \{ |t -s|^{\alpha H_t} , \ \ \ \big|{H_t  }- {H_s } \big|^{\alpha} \Big\} .\nonumber
\end{eqnarray}
Hence,
\begin{eqnarray*}
\Big|\Big| X_{H_t ,\alpha ,\lambda }( t) - X_{H_t ,\alpha ,\lambda }( v)\Big|\Big|_\alpha  \leq (2C_1) ^{1/\alpha}   \Big(  |t -s|^{  H_t}   + \big|{H_t  }- {H_s } \big|  \Big )  ,
\end{eqnarray*}
which gives  the desired upper bound in (\ref{f1d7}).

Next, consider the lower bound of $\big|\big| X_{H_t ,\alpha ,\lambda }( t) - X_{H_s ,\alpha ,\lambda }(s)\big|\big|_\alpha $.
 Write
\begin{eqnarray}
 \int_{-\infty}^{\infty}  \Big|  G_{H_t ,\alpha ,\lambda }(t, x)-G_{H_s ,\alpha ,\lambda }(s, x)   \Big|^{\alpha  } dx
&&\geq \int_{s}^{t } \Big|   e^{-\lambda(t-x) }(t-x) ^{H_t-\frac1{\alpha}}  \Big|^\alpha  dx   \nonumber \\
&&\geq e^{-\lambda \alpha (t-s) }\int_{s}^{t }(t-x) ^{H_t \alpha -1  }  dx   \nonumber \\
&&\geq e^{-\lambda \alpha (t-s) }\int_{s}^{t }(t-x) ^{H_t \alpha -1  }  dx  \nonumber \\
&&\geq e^{-\lambda \alpha }\frac1{H_t\alpha } (t-s) ^{H_t \alpha  }\nonumber \\
&&\geq e^{-\lambda \alpha }\frac1{b \alpha } (t-s) ^{H_t \alpha  } \nonumber
\end{eqnarray}
uniformly for all $t, s$ satisfying $s \leq t\leq s +1$.  Therefore, we have
\begin{eqnarray*}
\Big|\Big| X_{H_t ,\alpha ,\lambda }( t) - X_{H_s ,\alpha ,\lambda }( s)\Big|\Big|_\alpha   \geq   \,  e^{-\lambda} \Big(\frac1{b\alpha }\Big)^{1/\alpha}   |t -s|^{H_t }  ,
\end{eqnarray*}
which gives  the desired lower bound in (\ref{f1d7}).
  \hfill\qed

For $\alpha \in (0, 1],$ the next proposition shows that the upper bound of (\ref{f1d7}) is also exact.
 \begin{proposition} \label{hlmsg}
Assume $\alpha \in (0, 1]$ and $t_0>0.$  Then
there is a positive   number  $c$,  depending only on $ a, b, \lambda, t_0$ and $\alpha$,
such that
\begin{eqnarray}\label{f1d8}
\Big|\Big| X_{H_t ,\alpha ,\lambda }( t) - X_{H_s ,\alpha ,\lambda }(s)\Big|\Big|_\alpha  \geq c \, \Big(  |t -s|^{H_t}   + \big|{H_t  }- {H_s } \big|  \Big )
\end{eqnarray}
for all $t, s  $  satisfying $t_0\leq s \leq t \leq s+ 1$.
\end{proposition}
\emph{Proof.}
If $|t -s|^{H_t} \geq c_1 \big|{H_t  }- {H_s }\big|$ for some $c_1>0,$ depending only on $ a, b, \lambda, t_0$ and $\alpha$, then  (\ref{f1d7}) implies   (\ref{f1d8}). Otherwise,  we have
 \begin{eqnarray}  \label{tnlm}
 \big|{H_t  }- {H_s }\big| /|t -s|^{H_t}  \rightarrow \infty
 \end{eqnarray}
as $|t -s| \rightarrow 0.$
 Applying  the inequality
  $$\Big||x|^\alpha - |y|^\alpha \Big| \leq    |x-y|^\alpha, \ \ \ \ x, y \in \mathbf{R}\ \textrm{and} \ \alpha \in (0, 1], $$
 we have
\begin{eqnarray}
&&\int_{-\infty}^{\infty}  \Big|  G_{H_t ,\alpha ,\lambda }(t, x)-G_{H_s ,\alpha ,\lambda }(s, x)   \Big|^{\alpha  } dx \nonumber \\
&&\geq    \int_{0}^{s } \Big|   e^{-\lambda(t-x) }(t-x) ^{H_t-\frac1{\alpha}}  - e^{-\lambda(s-x) }(s-x) ^{H_s-\frac1{\alpha}}\Big|^\alpha  dx \nonumber \\
&&=   \int_{0}^{s } \Big|  e^{-\lambda(t-x) }(t-x) ^{H_t-\frac1{\alpha}} - e^{-\lambda(s-x) }(s-x) ^{H_t-\frac1{\alpha}} \nonumber \\
&&  \ \ \ \ \ \  \ \ \ \ \ \  + \,  e^{-\lambda(s-x) }(s-x) ^{H_t-\frac1{\alpha}}  - e^{-\lambda(s-x) }(s-x) ^{H_s-\frac1{\alpha}}\Big|^\alpha  dx \nonumber \\
&&\geq    \int_{0}^{s } \Big|  e^{-\lambda(s-x) }(s-x) ^{H_t-\frac1{\alpha}}  - e^{-\lambda(s-x) }(s-x) ^{H_s-\frac1{\alpha}}\Big|^\alpha   dx \nonumber \\ && \ \ \ \ \ \  \ \ \ \ -\int_{0}^{s } \Big|  e^{-\lambda(t-x) }(t-x) ^{H_t-\frac1{\alpha}} - e^{-\lambda(s-x) }(s-x) ^{H_t-\frac1{\alpha}}\Big|^\alpha   dx. \nonumber
\end{eqnarray}
By the mean value theorem and (\ref{ineqs12}), the last inequality implies that   for $\alpha \in (0, 1]$ and all $t, s$ satisfying $0\leq s \leq t \leq s+1$,
\begin{eqnarray*}
&&\int_{-\infty}^{\infty}  \Big|  G_{H_t ,\alpha ,\lambda }(t, x)-G_{H_s ,\alpha ,\lambda }(s, x)   \Big|^{\alpha  } dx \nonumber \\
&& \geq  \int_{0}^{s} e^{-\lambda \alpha (s-x) }  (s-x) ^{ \alpha H_\theta -1} |\log(s-x)|^{\alpha}  dx \, \Big|H_t- H_s\Big|^\alpha       - C_{11}    |t-s|^{\alpha H_t}       \nonumber  \\
&& \geq  c_{00}  \Big|H_t- H_s\Big|^\alpha       - C_{11}    |t-s|^{\alpha H_t},
\end{eqnarray*}
where $c_{00}, C_{11} >0$ depending only on $ a, b, \lambda, t_0$ and $\alpha$. By (\ref{tnlm}),  it follows  that   for $\alpha \in (0, 1],$
\begin{eqnarray*}
 \int_{-\infty}^{\infty}  \Big|  G_{H_t ,\alpha ,\lambda }(t, x)-G_{H_s ,\alpha ,\lambda }(s, x)   \Big|^{\alpha  } dx
  \geq c   \Big|H_t- H_s\Big|^\alpha ,
\end{eqnarray*}
where $c >0$ depending only on $ a, b, \lambda, t_0$ and $\alpha$.
Therefore (\ref{f1d8}) holds.
  \hfill\qed

\section{Localisability  and strong localisability } \label{sec8}

Recall that a stochastic process $X(t), t \in  \mathbf{R},$ is said to be $h-$localisable at $u$ (cf.\ Falconer \cite{F,F2}),
with $h>0$, if there exists a non-trivial process $X_u '$, called the tangent process of $X$ at $u$, such that
\begin{eqnarray}\label{indis}
\lim_{r\searrow 0} \frac{X(u+rv)-X(u)}{r^h} \stackrel{fdd}{= } X_u'(v),
\end{eqnarray}
where  $\stackrel{fdd}{= }$ stands for convergence in  finite-dimensional distributions.

The following proposition shows that LTFmSM is $H-$localisable.
\begin{proposition}\label{pr2.1} Assume that $\alpha(x)$ is continuous on $\mathbf{R}$. When $ 1 /a <  H  <1,$ the LTFmSM process $X_{H ,\alpha(x),\lambda }( t) $ is $H-$localisable at $u  $   with local form
\begin{eqnarray}\label{ffgls}
  X (t):=\int_{-\infty}^{\infty}\Big[  (t-x)_+^{H -\frac1{\alpha(u)}} -  (-x)_+^{H-\frac1{\alpha(u)}}\Big] dZ_{\alpha(u)}(x),
\end{eqnarray}
where $dZ_{\alpha(u)}(x)$ is a symmetric $\alpha(u)-$stable random measure.
\end{proposition}
\emph{Proof.}
 Given $ u_1 < u_2 < ...< u_d,$   denote
\begin{eqnarray}
S_r(u_k)= \frac{X_{H,\alpha ,\lambda }( u + ru_k) - X_{H ,\alpha ,\lambda }( u )}{ r^{H }} \nonumber
\end{eqnarray}
for $r>0$ and  $k=1,...,d.$ Then
\begin{eqnarray*}
 &&\mathbf{E}\Big[e^{i\sum_{k=1}^d \theta_k S_r(u_k)} \Big] \\
 &&=\  \exp\Bigg\{ - \int_{-\infty}^{\infty} \Big|\sum_{k=1}^d\theta_kr^{-H} \Big(G_{H ,\alpha(x),\lambda }(u + ru_k, x)-G_{H ,\alpha(x),\lambda }(u, x) \Big) \Big|^{\alpha(x)} dx\Bigg\}.
\end{eqnarray*}
Let $x=u + rz.$ It follows that
\begin{eqnarray*}
&& \int_{-\infty}^{\infty} \Big|\sum_{k=1}^d\theta_kr^{-H} \Big(G_{H ,\alpha(x),\lambda }(u + ru_k, x)-G_{H ,\alpha(x),\lambda }(u, x) \Big) \Big|^{\alpha(x)} dx \\ &&=   \int_{-\infty}^{\infty} \Big|  \sum_{k=1}^d\theta_k  \Big(e^{-\lambda r( u_k-z)_+}(  u_k-z)_+^{H -\frac1{\alpha(u + rz)}}  -e^{-\lambda r (   -z)_+}(  -z)_+^{H -\frac1{\alpha(u+rz)}} \Big) \Big|^{\alpha(u+rz)} dz.
\end{eqnarray*}
Recall that $\alpha(x)  \in [a, b] $ is a continuous function  on $\mathbf{R}.$
Thus
\begin{eqnarray*}
  && \lim_{r\rightarrow 0}\Big|  \sum_{k=1}^d\theta_k  \Big(e^{-\lambda r( u_k-z)_+}(  u_k-z)_+^{H -\frac1{\alpha(u + rz)}}  -e^{-\lambda r (   -z)_+}(  -z)_+^{H -\frac1{\alpha(u+rz)}} \Big) \Big|^{\alpha(u+rz)} \\
  && = \Big|  \sum_{k=1}^d \theta_k  \Big( (  u_k-z)_+^{H -\frac1{\alpha(u  )}}  - (  -z)_+^{H -\frac1{\alpha(u )}} \Big) \Big|^{\alpha(u )}.
\end{eqnarray*}
It is obvious that for $z<\min\{u_1, 0\}-1,$ $0<r<1$ and  $1/a <  H < 1$,
\begin{eqnarray*}
  && \Big|  \sum_{k=1}^d\theta_k  \Big(e^{-\lambda r( u_k-z)_+}(  u_k-z)_+^{H -\frac1{\alpha(u + rz)}}  -e^{-\lambda r (   -z)_+}(  -z)_+^{H -\frac1{\alpha(u+rz)}} \Big) \Big|^{\alpha(u+rz)} \\
  && \leq   \sum_{k=1}^d |\theta_k|^{\alpha(u+rz)}  \Big| e^{-\lambda r  u_k } (  u_k-z)^{H - \frac 1{\alpha(u + rz)}} - (  -z)^{H - \frac{1}{\alpha(u + rz)}}  \Big|^{\alpha(u+rz)}  \, \, \,  \\
  && \leq   \sum_{k=1}^d |\theta_k|^{\alpha(u+rz)}  e^{ \lambda r  |u_k| } \Big|  (  u_k-z)^{H - \frac 1{\alpha(u + rz)}} - (  -z)^{H - \frac{1}{\alpha(u + rz)}}  \Big|^{\alpha(u+rz)}     \\
  && \leq   \sum_{k=1}^d \Big( |\theta_k|^{a } + |\theta_k|^{b }\Big) e^{ \lambda   |u_k| } \sup_{\alpha \in [a, b]} \Big| (  u_k-z)^{H - 1/\alpha } - (  -z)^{H - 1/\alpha }  \Big|^\alpha  \\
  && \leq   \sum_{k=1}^d \Big( |\theta_k|^{a } + |\theta_k|^{b }\Big) e^{ \lambda   |u_k| } \Big|H-\frac1a \Big|^a  \Big( \min\{u_1, 0\}-z \Big)^{H a-1-a    } (|u_k|^a + |u_k|^b),
\end{eqnarray*}
and that for $ z\geq \min\{u_1, 0\}-1$ and $0<r<1,$
\begin{eqnarray*}
  && \Big|  \sum_{k=1}^d\theta_k  \Big(e^{-\lambda r( u_k-z)_+}(  u_k-z)_+^{H -\frac1{\alpha(u + rz)}}  -e^{-\lambda r (   -z)_+}(  -z)_+^{H -\frac1{\alpha(u+rz)}} \Big) \Big|^{\alpha(u+rz)} \\
  && \leq   \sum_{k=1}^d |\theta_k|^{\alpha(u+rz)}  \Big( e^{-\lambda r   (  u_k-z)_+ \alpha(u+rz)} (  u_k-z)_+^{H \alpha(u + rz)- 1}  + e^{-\lambda r (   -z)_+\alpha(u + rz)} (  -z)_+^{H \alpha(u + rz)- 1}  \Big)    \\
   && \leq   \sum_{k=1}^d |\theta_k|^{\alpha(u+rz)}  \Big(   (  u_k-z)_+^{H \alpha(u + rz)- 1}  +  (  -z)_+^{H \alpha(u + rz)- 1}  \Big)    \\
  && \leq   \sum_{k=1}^d \Big( |\theta_k|^{a } + |\theta_k|^{b }\Big)\Big(   (  u_k-z)_+^{H a- 1} + (  u_k-z)_+^{H b- 1} +  (  -z)_+^{H a- 1}+ (  -z)_+^{H b- 1}  \Big).
\end{eqnarray*}
The dominated convergence theorem implies that
\begin{eqnarray*}
&& \lim_{r\rightarrow 0}\mathbf{E}\Big[e^{i\sum_{k=1}^d \theta_k S_r(u_k)} \Big]  \\
  && =\  \exp\Bigg\{ - \int_{-\infty}^{\infty}\Big|  \sum_{k=1}^d \theta_k  \Big( (  u_k-z)_+^{H -\frac1{\alpha(u  )}}  - (  -z)_+^{H -\frac1{\alpha(u )}} \Big) \Big|^{\alpha(u )} dz\Bigg\}\\
  && = \mathbf{E}\Big[e^{i\sum_{k=1}^d \theta_k X (u_k)} \Big] ,
\end{eqnarray*}
where $X (\cdot)$ is defined  by (\ref{ffgls}).
By L\'{e}vy's continuous theorem, we have
 \begin{eqnarray*}
 \lim_{r\rightarrow 0} S_r(u_k)\stackrel{fdd}{=}   X (u_k)  .
\end{eqnarray*}
Thus   $  X_{H ,\alpha(x),\lambda }( t), t \in \mathbf{R},$ is $H$-localisable at $u  $  to $X (\cdot)$ defined  by (\ref{ffgls}).  \hfill\qed

Recall that $X(t), t \in  \mathbf{R},$ is said to be $h$-strongly localisable at $u$
to $X_u'(v)$ with $h>0$  (cf.\ Falconer and  Liu \cite{FL12}),
if the convergence in (\ref{indis}) occurs in distribution with respect to the metric of
uniform convergence on bounded intervals, and  $X$ and $X_t'$ have versions in $C(R)$ (the space of continuous function
on $\mathbf{R}$).

The next proposition shows that when $1/a < H < 1,$ the LTFmSM is $H$-strongly localisable.
\begin{proposition}\label{th3.7}
Assume that $\alpha(x)$ is continuous on $\mathbf{R}$. When $ 1 /a <  H  <1,$    the process $ X_{H ,\alpha(x),\lambda }( t) $ is $H-$strongly localisable at $u$ to the LFSM defined  by (\ref{ffgls}).
\end{proposition}
\emph{Proof.}
By Theorem 3.2 of   Falconer  and  Liu \cite{FL12}, it is sufficient to prove that for each bounded interval $J,$ there is a positive $r_0$ such that for any $r \in (0, r_0),$ $$ \int_{-\infty}^{\infty}  \Big| \frac{G_{H ,\alpha(x),\lambda }(u+rt, x)-G_{H ,\alpha(x),\lambda }(u+rv, x)}{r^H} \Big|^{\alpha(x)}dx \leq C\, |t-v|^{aH}, \ \ \   \ \ \ t,v \in J,$$
where $C$ is a constant. Indeed, by (\ref{sjkns}), for any $0< r \leq \min\{ 1/|t-v|,  1\},$ we have
\begin{eqnarray*}
&&\int_{-\infty}^{\infty}  \Big| \frac{G_{H ,\alpha(x),\lambda }(u+rt, x)-G_{H ,\alpha(x),\lambda }(u+rv, x)}{r^H} \Big|^{\alpha(x)}dx \\
&& \leq \frac{1}{r^{Ha}}\int_{-\infty}^{\infty}  \Big|G_{H ,\alpha(x),\lambda }(u+rt, x)-G_{H ,\alpha(x),\lambda }(u+rv, x) \Big|^{\alpha(x)}dx \\
&& \leq \frac{1}{r^{Ha}} C\, |rt-rv|^{Ha} = C\, |t-v|^{aH}.
\end{eqnarray*}
This completes the proof of Proposition \ref{th3.7}. \hfill\qed

When $\lambda=0$ and $ 1/a  < H < 1+ 1/b -1/a,$ Falconer  and  Liu  proved that  $  X_{H ,\alpha(x), 0 }( t), t \in \mathbf{R}, $ is $H-$strongly localisable, see Proposition 4.3 of \cite{FL12}. Now Proposition \ref{th3.7} extends the result  of Falconer  and  Liu.

\section*{Acknowledgments}
We would like to thank the anonymous referees and the editor for their helpful comments and suggestions.
Jacques L\'{e}vy V\'{e}hel gratefully acknowledges support from SMABTP.  The work of Xiequan Fan has been partially supported by the
 National Natural Science Foundation of China (Grant nos.\ 11601375 and 11626250).


\end{document}